\documentclass[11pt,twoside]{amsart}
\usepackage{}
\usepackage{amsthm}
\usepackage{amsmath}
\usepackage{amssymb}
\usepackage{amsfonts}
\usepackage{latexsym}
\usepackage{enumitem}
\usepackage{mathrsfs}
\usepackage[all]{xy} \SelectTips{eu}{}
\usepackage{hyperref}
\usepackage{eucal}
\usepackage{xcolor}
\usepackage{tikz,tikz-cd}
\usetikzlibrary{arrows}
\usetikzlibrary{matrix}

\newtheorem{intthm}{Theorem}[]

\newtheorem*{intque*}{Question}
\newtheorem*{intexa*}{Example}
\newcommand{\numberseries}{\bfseries}   
\newlength{\thmtopspace}                
\newlength{\thmbotspace}                
\newlength{\thmheadspace}               
\newlength{\thmindent}                  
\newtheoremstyle{bfupright head,slanted body}
{\thmtopspace}{\thmbotspace}
{\slshape}{\thmindent}{\bfseries}{.}{\thmheadspace}
{{\numberseries \thmnumber{#2\;}}\thmnote{#3}}

\newtheoremstyle{bfupright head,upright body}
{\thmtopspace}{\thmbotspace}
{\upshape}{\thmindent}{\bfseries}{.}{\thmheadspace}
{{\numberseries \thmnumber{#2\;}}\thmnote{#3}}

\newtheoremstyle{fixed bf head,slanted body}
{\thmtopspace}{\thmbotspace}{\slshape}
{\thmindent}{\bfseries}{.}{\thmheadspace}
{{\numberseries \thmnumber{#2\;}}\thmname{#1}\thmnote{ (#3)}}

\newtheoremstyle{fixed bf head,upright body}
{\thmtopspace}{\thmbotspace}{\upshape}
{\thmindent}{\bfseries}{.}{\thmheadspace}
{{\numberseries \thmnumber{#2\;}}\thmname{#1}\thmnote{ (#3)}}

\newtheoremstyle{numbered paragraph}
{\thmtopspace}{\thmbotspace}{\upshape}
{\thmindent}{\upshape}{}{\thmheadspace}
{{\numberseries \thmnumber{#2.}}}

\theoremstyle{bfupright head,slanted body}
\newtheorem{res}{}[section]             \newtheorem*{res*}{}

\theoremstyle{bfupright head,upright body}
\newtheorem{bfhpg}[res]{}               \newtheorem*{bfhpg*}{}

\theoremstyle{fixed bf head,slanted body}

\newtheorem*{Apthm*}{Approximation Theorem}
         \newtheorem*{thm*}{Theorem}
\newtheorem{prp}[res]{Proposition}      \newtheorem*{prp*}{Proposition}
        \newtheorem*{cor*}{Corollary}
\newtheorem{lem}[res]{Lemma}            \newtheorem*{lem*}{Lemma}
\newtheorem{que}[res]{Question}         \newtheorem*{que*}{Question 1}
\theoremstyle{fixed bf head,upright body}
\newtheorem{dfn}[res]{Definition}       \newtheorem*{dfn*}{Definition}
\newtheorem{rmk}[res]{Remark}           \newtheorem*{rmk*}{Remark}
\newtheorem{exa}[res]{Example}           \newtheorem*{exa*}{Example}
           \newtheorem*{nota*}{Notation}
           \newtheorem*{setup*}{Setup}
\newtheorem{setup and notation}[res]{Setup and notation}  \newtheorem*{setup and notation*}{Setup and notation}

\theoremstyle{numbered paragraph}
\newtheorem{ipg}[res]{}

\newlength{\thmlistleft}        
\newlength{\thmlistright}       
\newlength{\thmlistpartopsep}   
\newlength{\thmlisttopsep}      
\newlength{\thmlistparsep}      
\newlength{\thmlistitemsep}     

\newcounter{eqc}
	{\end{list}}%

\newcounter{prt}
\newenvironment{prt}{\begin{list}{\upshape (\alph{prt})}%
		{\usecounter{prt}%
			\setlength{\leftmargin}{2.5em}
			\setlength{\labelwidth}{\thmlistleft}%
			\setlength{\rightmargin}{\thmlistright}%
			\setlength{\partopsep}{\thmlistpartopsep}%
			\setlength{\topsep}{\thmlisttopsep}%
			\setlength{\parsep}{\thmlistparsep}%
			\setlength{\itemsep}{\thmlistitemsep}}}%
	{\end{list}}%

\newcounter{rqm}
\newenvironment{rqm}{\begin{list}{\upshape (\arabic{rqm})}%
		{\usecounter{rqm}%
			\setlength{\leftmargin}{2.5em}
			\setlength{\labelwidth}{\thmlistleft}%
			\setlength{\rightmargin}{\thmlistright}%
			\setlength{\partopsep}{\thmlistpartopsep}%
			\setlength{\topsep}{\thmlisttopsep}%
			\setlength{\parsep}{\thmlistparsep}%
			\setlength{\itemsep}{\thmlistitemsep}}}%
	{\end{list}}%

	{\end{list}}%

\newenvironment{prf*}[1][Proof]{%
	\begin{proof}[\bf #1]
		\setcounter{equation}{0}
		}
	{\end{proof}
}

\newcommand{\pgref}[1]{\ref{#1}}

\renewcommand{\eqref}[1]{(\pgref{eq:#1})}

\newcommand{\pagecite}[2][?]{\cite[page.~#1]{#2}}
\newcommand{\thmcite}[2][?]{\cite[Theorem~#1]{#2}}
\newcommand{\prpcite}[2][?]{\cite[Proposition~#1]{#2}}

\newcommand{\corcite}[2][?]{\cite[Corollary~#1]{#2}}
\newcommand{\lemcite}[2][?]{\cite[Lemma~#1]{#2}}

\numberwithin{equation}{res}

\def\urltilda{\kern -.15em\lower .7ex\hbox{\~{}}\kern .04em}

\newcommand{\Ker}[1]{\nobreak{\operatorname{ker}#1}}
\newcommand{\Coker}[1]{\nobreak{\operatorname{coker}#1}}
\newcommand{\Mon}{\mathsf{Mon}}
\newcommand{\Epi}{\mathsf{Epi}}
\newcommand{\sfE}{\mathsf{E}}
\newcommand{\sfM}{\mathsf{M}}
\newcommand{\ps}{^\mathsf{ps}}

\newcommand{\Cat}{\mathsf{Cat}}
\newcommand{\MCat}{\mathsf{ModCat}}

\newcommand{\sfA}{\mathsf{A}}
\newcommand{\sfC}{\mathsf{C}}
\newcommand{\sfF}{\mathsf{F}}
\newcommand{\sfL}{\mathsf{L}}
\newcommand{\sfW}{\mathsf{W}}
\newcommand{\sfX}{\mathsf{X}}
\newcommand{\sfCof}{\mathsf{Cof}}
\newcommand{\sfInj}{\mathsf{Inj}}
\newcommand{\Rep}[2]{{\sf Rep}(#1,#2)}
\newcommand{\sfR}{_\mathsf{R}}

\newcommand{\scrF}{\mathscr{F}}
\newcommand{\scrR}{\mathscr{R}}

\newcommand{\calB}{\mathcal{B}}

\newcommand{\calM}{\mathcal{M}}
\newcommand{\calI}{\mathcal{I}}
\newcommand{\calT}{\mathcal{T}}

\newcommand{\id}{\mathrm{id}}

\newcommand{\Fun}{\mathsf{Fun}}

\newcommand{\WFib}{\mathsf{Wald.opFib}}

\newcommand{\Wald}{\mathsf{WaldCat}}

\newcommand{\opfib}{\mathsf{opFib}}

\newcommand{\Mor}{\mathsf{Mor}}
\newcommand{\cMor}{\mathsf{coMor}}
\newcommand{\coE}{\mathsf{coE}}
\newcommand{\Ob}{\mathsf{Ob}}

\def\soft#1{\leavevmode\setbox0=\hbox{h}\dimen7=\ht0\advance
	\dimen7 by-1ex\relax\if t#1\relax\rlap{\raise.6\dimen7
		\hbox{\kern.3ex\char'47}}#1\relax\else\if T#1\relax
	\rlap{\raise.5\dimen7\hbox{\kern1.3ex\char'47}}#1\relax
	\else\if d#1\relax\rlap{\raise.5\dimen7\hbox{\kern.9ex
			\char'47}}#1\relax\else\if D#1\relax\rlap{\raise.5\dimen7
		\hbox{\kern1.4ex\char'47}}#1\relax\else\if l#1\relax
	\rlap{\raise.5\dimen7\hbox{\kern.4ex\char'47}}#1\relax
	\else\if L#1\relax\rlap{\raise.5\dimen7\hbox{\kern.7ex
			\char'47}}#1\relax\else\message{accent \string\soft
		\space #1 not defined!}#1\relax\fi\fi\fi\fi\fi\fi}
\makeatletter
\def\part{\@startsection{part}{1}%
\z@{.7\linespacing\@plus\linespacing}{.8\linespacing}%
{\LARGE\sffamily\centering}}

\makeatletter
\def\l@section{\@tocline{1}{2pt}{0pc}{}{}}
\makeatother
\let\oldtocpart=\tocpart
\renewcommand{\tocpart}[2]{\bf\large\oldtocpart{#1}{#2}}
\let\oldtocsection=\tocsection
\renewcommand{\tocsection}[2]{\bf\oldtocsection{#1}{#2}}

\title[Constructions of Waldhausen categories via Grothendieck opfibrations]
{Constructions of Waldhausen categories via Grothendieck opfibrations}

\date{ \today}

\keywords{Waldhausen category (structure); cofibration;
weak equivalence; Grothendieck opfibration; left rooted quiver.}

\subjclass[2010]{18G25; 18A25; 18A40}

\thanks{Z.X. Di was partly supported by NSF of China (Grant No. 11971388),
the Scientific Research Funds of Huaqiao University (Grant No. 605-50Y22050) and
the Fujian Alliance of Mathematics (Grant No. 2024SXLMMS04);
L.P. Li was partly supported by NSF of China (Grant No. 12171146);
L. Liang was partly supported by NSF of China (Grant No. 12271230).}

\author[Z.X. Di]{Zhenxing Di}
\address{Z.X. Di \ School of Mathematical Sciences, Huaqiao University, Quanzhou 362021, China}
\email{dizhenxing@163.com}

\author[L.P. Li]{Liping Li}
\address{L.P. Li \ LCSM (Ministry of Education), Department of Mathematics, Hunan Normal University, Changsha 410081, China.}
\email{lipingli@hunnu.edu.cn}

\author[L. Liang]{Li Liang}
\address{L. Liang \ Department of Mathematics, Lanzhou Jiaotong University, Lanzhou 730070, China}
\email{lliangnju@gmail.com}
\urladdr{https://sites.google.com/site/lliangnju}

\begin{document}

\begin{abstract}
Given a Grothendieck opfibration $p: \calT \to \calB$,
we describe a method to construct a Waldhausen category structure on the total category $\calT$ via combining Waldhausen category structures on the fibers $\calT_A$ for $A \in \Ob(\calB)$ and the basis category $\calB$. As an application, we show that if $\sfE$ is a Waldhausen category with small coproducts such that the class of cofibrations is the left part of a weak factorization system in $\sfE$, then the representation category $\Rep{Q}{\coE}$ of a left rooted quiver $Q$ is a Waldhausen category, where $\coE$ is the subcategory of $\sfE$ whose morphisms are cofibrations.
\end{abstract}

\maketitle

\section{Introduction}\label{intro}
\noindent
Given a pseudo functor $\scrF$ from a category $\calB$ to the (meta) 2-category $\Cat$ of categories, the lax colimit of $\scrF$ is represented by its \emph{Grothendieck construction} $\int_\calB \scrF$. The association of $\scrF$ to the canonical projection
$p : \int_\calB \scrF \to \calB$ exhibits one of the most fundamental relations in category theory, namely, the following equivalence between pseudo functors from $\calB$ to $\Cat$ and Grothendieck opfibrations over $\calB$:
\begin{equation}\label{groth corre}
\begin{gathered}
\tag{\ref{intro}.1}
\int : \Fun\ps(\calB, \Cat) \to \opfib(\calB),
\end{gathered}
\end{equation}
which is now known as \emph{Grothendieck correspondence} \cite{Groth61}.

A natural question one may consider is to establish similar correspondences for categories with extra structure. This question has been extensively studied in the literature, see for instances \cite{MoVa20, BeWo19, Lanf23} for monoidal categories, enriched categories and tangent categories. For model categories, Harpaz and Prasma \cite{HP2015} showed that a covariant pseudo functor satisfying certain conditions from a model category $\calM$ to $\MCat$, the category of model categories and Quillen adjunctions, induces the integral model structure on its Grothendieck construction, and established an equivalence between these pseudo functors and model bifibrations over $\calM$. Inspired by \cite{HP2015} as well as Roig \cite{Roig94} and Stanculescu \cite{Stan12}, Cagne and Melli\`es identified in \cite{CM20} necessary and sufficient conditions on a Grothendieck bifibration $p : \calT \to \calB$ to ensure that model structures on the fibers $\calT_A$ for $A \in \Ob(\calB)$ and the basis category $\calB$ combine into a model structure on the total category $\calT$. Furthermore, via this Grothendieck construction, they revisited the traditional definition of Reedy model structures and exhibited their bifibrational nature.

\emph{Waldhausen categories} were introduced in \cite{Wald83} to extend higher algebraic $K$-theory brought by Quillen \cite{QuiK72}. Note that Waldhausen categories only involve cofibrations and weak equivalences, and hence in general are different from model categories. Thus one may ask the following question:

\begin{que}\label{que 1}
Given a pseudo functor $\scrF : \calB \to \Wald$, where $\calB$ is a Waldhausen category and $\Wald$ is the category of Waldhausen categories and exact functors, how to endow a Waldhausen category structure on its Grothendieck construction $\int_\calB \scrF$?
\end{que}

Let $\WFib (\calB)$ be the subcategory of $\opfib(\calB)$ consisting of Grothendieck opfibrations $p : \calT \to \calB$ satisfying the following conditions:
\begin{itemize}
\item the fiber $\calT_A$ is a Waldhausen category
      for each object $A$ in $\calB$;

\item the reindexing functor $u_!: \calT_A\to \calT_B$ is exact
      for each morphism $u$ in $\calB$.
\end{itemize}
One may refer to objects in $\WFib (\calB)$ as \emph{Waldhausen opfibrations}. It is clear that (\ref{groth corre}) restricts to an equivalence
\[
\int_\calB : \Fun\ps(\calB, \Wald) \to \WFib(\calB),
\]
where $\Fun\ps(\calB, \Wald)$ is the subcategory of $\Fun\ps(\calB, \Cat)$ spanned by pseudo functors whose domain lies in $\Wald$. This equivalence allows us to answer the above question from the point of view of Waldhausen opfibrations, which is also beneficial to giving further applications in representation categories of left rooted quivers.

We would like to show that for a Waldhausen opfibration $p : \calT \to \calB$ with $\calB$ a Waldhausen category, the total category $\calT$ inherits a Waldhausen category structure from the ones assigned to the basis category $\calB$ and to the fibers $\calT_A$ for each object $A$ in $\calB$. Note that any morphism $f : X \to Y$ in $\calT$ that is above a morphism
$u : A \to B$ in $\calB$ can be factored uniquely as a cocartesian lifting of $X$ along $u$ followed by a morphism $f_\triangleright$ in the fiber $\calT_B$ by the cocartesian universal property of the lifting (see \ref{factorization 1}). One can always equip $\calT$ with two classes of total cofibrations and total weak equivalences described as:
\begin{itemize}
\item
A \emph{total cofibration} (resp., \emph{total weak equivalence}) is a morphism
$f : X \to Y$ in $\calT$ that is above a cofibration (resp., \emph{weak equivalence})
$u : A \to B$ in $\calB$ such that $f_\triangleright$ is a cofibration
(resp., \emph{weak equivalence}) in the fiber $\calT_B$.
\end{itemize}
The following is the first main result of the paper.

\begin{intthm}\label{total E IS WALD}
Let $p: \calT \to \calB$ be a Waldhausen opfibration. If the basis category $\calB$ is a Waldhausen category, then the classes of total cofibrations and total weak equivalences given above endow the total category $\calT$ with a Waldhausen category structure.
\end{intthm}

Since the canonical projection $p : \int_\calB \scrF \to \calB$ associated to the pseudo functor $\scrF$ in Question \ref{que 1} is clearly a Waldhausen opfibration, Theorem \ref{total E IS WALD} implies that the Grothendieck construction $\int_\calB \scrF$ considered in Question \ref{que 1} forms a Waldhausen category, where a morphism $(u, \phi) : (A, X) \to (B, Y)$ is a cofibration (resp., weak equivalence) if $u : A \to B$ is a cofibration (resp., weak equivalence) in $\calB$ and
$\phi : \scrF_u(X) \to Y$ is a cofibration (resp., weak equivalence) in $\scrF_B$.

Let $\sfE$ be a Waldhausen category and $\coE$ its subcategory with the same objects whose morphisms are cofibrations. Waldhausen showed in \cite{Wald83} that both the morphism category $\Mor(\sfE)$ and its full subcategory $\cMor(\sfE)$ consisting of cofibrations form also Waldhausen categories; see Example \ref{classical exa mor}. As illustrations of Theorem \ref{total E IS WALD}, we reobtain the Waldhausen category structures on $\Mor(\sfE)$ and $\cMor(\sfE)$ via the codomain and domain functors, respectively; see Examples \ref{reobtain classical exa1} and \ref{rebtain classical 2}. Since $\Mor(\sfE)$ (resp., $\cMor(\sfE)$) can be identified with the functor category $\Fun(\calI,\sfE)$ (resp., $\Fun(\calI, \coE)$), where $\calI$ is the free category associated to the quiver with two vertexes and one arrow, one may ask the following question:
\begin{que}\label{que 2}
Let $\calI$ be a small index category and $\sfE$ a Waldhausen category. Is there a canonical Waldhausen structure on $\Fun(\calI, \sfE)$ or $\Fun(\calI, \coE)$ generalizing the above two examples?
\end{que}

The answer for the first functor category is transparent. Indeed, define a morphism $f : X \to Y$ in $\Fun(\calI,\sfE)$ to be a cofibration (resp., weak equivalence) if $f_i : X_i \to Y_i$ is a cofibration (resp., weak equivalence) in $\sfE$ for each object $i$ in $\calI$. Since isomorphisms and pushouts in $\Fun(\calI,\sfE)$ are determined componentwise, it is easy to check that $\Fun(\calI,\sfE)$ forms a Waldhausen category in this way. However, the situation for the second functor category becomes much more complicated, and at this moment we cannot figure out a completely satisfactory solution. Instead, in this paper we will focus on the special case that $\calI$ is the free category associated to a left rooted quiver $Q$ (see Section \ref{Appl in rep cat} for the definition of left rooted quivers). Using Theorem \ref{total E IS WALD} and the approach of Cagne and Melli\`es \cite{CM20} for Reedy model structures, we show that if $\sfE$ has small coproducts and the class of cofibrations is the left part of a weak factorization system in $\sfE$, then $\Rep{Q}{\coE} = \Fun(\calI, \coE)$ admits a Waldhuasen category structure.

Before displaying the explicit statement of the result,
let us describe the key steps and necessary notation.
As one might expect, the proof of the result goes through a transfinite induction and
the most technical part of the proof is the induction step for successor ordinals.
Let $\{V_{\mu}\}_{\mu \leqslant \zeta}$ be the transfinite sequence of subsets of vertexes in $Q$, and let $Q_{\mu} = (V_{\mu}, \Gamma_{\mu})$ denote the subquiver of $Q$ spanned by $V_{\mu}$ for each ordinal $\mu \leqslant \zeta$.
It is shown that the restriction functor
$\iota_\mu^* : \Rep{Q_{\mu+1}}{\coE} \to \Rep{Q_{\mu}}{\coE}$
with respect to the inclusion $\iota_\mu : Q_\mu \to Q_{\mu+1}$
is a Waldhausen opfibration (see Proposition \ref{iota mu is a Waldop})
under which the induction step for a successor ordinal
can be completed by Theorem \ref{total E IS WALD};
see \ref{pf of thB} for details.
The key role playing in the proof of the above result is that
the fiber $\Rep{Q_{\mu+1}}{\coE}_A$ of a representation $A \in \Rep{Q_{\mu}}{\coE}$
is isomorphic to the product category
$\prod_{i \in V_{\mu+1} \backslash V_{\mu}} \overline{\sfL_i(A)/\sfE}$
(see Proposition \ref{isomorphism of fibe}),
where $\sfL_i(A) = \oplus_{\alpha \in \Gamma_{\mu + 1}(\bullet, i)} A_{s(\alpha)}$
is indeed the latching object of $A$ at $i$ and
$\overline{\sfL_i(A)/\sfE}$ is the full subcategory of the undercategory $\sfL_i(A)/\sfE$
consisting of all cofibrations $(\sfL_i(A) \rightarrowtail \bullet)$ in $\sfE$,
which is a Waldhausen category via a natural way (see Example \ref{over cate is}(2)).
For any representation $X$ in $\Rep{Q}{\coE}$ and any vertex $i$ in $Q$,
by the universal property of the coproduct,
we obtain a morphism $\sfL_i(X) \to X_i$, which is a cofibration in $\sfE$;
one refers to \ref{1-1 for exten} for a similar proof.
Hence, for any morphism $f : X \to Y$ in $\Rep{Q}{\coE}$,
by the axiom (C3) in Definition \ref{df of wald cat and fun},
we get a pushout appearing as the inner square of the commutative diagram
\[
\xymatrix@R=1cm@C=1cm{
   \sfL_i(X)
     \ar@{>->}[d]_{}
     \ar[r]^-{\sfL_i(f)}
     \ar@{}[rd]^(0.6)>>{\lrcorner}
& \sfL_i(Y)
     \ar@{>->}[d]_{}
     \ar@/^0.8pc/[ddr]^{}                          \\
  X_i
     \ar[r]^-{}
     \ar@/_0.8pc/[drr]_{f_i}
&  \sfL_i(Y) \sqcup_{\sfL_i(X)} X_i
     \ar@{.>}[dr]|-{}                              \\
&& Y_i}
\]
of solid arrows in $\sfE$.
By the universal property of the pushout,
there is a natural morphism $\sfL_i(Y) \sqcup_{\sfL_i(X)} X_i \to Y_i$.
Now, we can state the second main result of the paper.

\begin{intthm}\label{MAIN FOR REP CAT}
Let $Q$ be a left rooted quiver and $\sfE$ a Waldhausen category with small coproducts. Denote by $\sfC$ $($resp., $\sfW$$)$ be the class of all cofibrations $($resp., weak equivalences$)$ in $\sfE$. If $(\sfC, \sfC^\Box)$ forms a weak factorization system in $\sfE$, then $\Rep{Q}{\coE}$ forms a Waldhausen category,
where a morphism $f: X \to Y $ is a cofibration $($resp., weak equivalence$)$ if the induced morphism
\[
X_i \sqcup_{\sfL_i(X)} \sfL_i(Y) \to Y_i
\]
lies in $\sfC$ $($resp., $\sfW$$)$ for each vertex $i$ in $Q$.
\end{intthm}

\begin{rmk}
If further $\sfW$ is closed under small coproducts and pushouts along morphisms in $\sfC$, and satisfies the 2-out-of-3 property for composable morphisms in $\sfE$,
then a morphism $f : X \to Y$ is a weak equivalence in $\Rep{Q}{\coE}$ if and only if $f_i : X_i \to Y_i$ lies in $\sfW$ for each vertex $i$ in $Q$ (see Remark \ref{pointed w}).
\end{rmk}

Given a Quillen model category $\sfM$ with small coproducts and $(\sfC, \sfW, \sfF)$ its model structure, we know that $\Rep{Q}{\sfM}$ forms also a Quillen model category;
see for instance Hirschhorn \cite{Hir03}.
By Example \ref{exa model} and Theorem \ref{MAIN FOR REP CAT},
both the full subcategories $\sfCof\sfR$ and $\Rep{Q}{\overline{\sfCof}}$
of $\Rep{Q}{\sfM}$ form Waldhausen categories,
where $\sfCof\sfR$ denotes the full subcategory of all cofibrant objects in $\Rep{Q}{\sfM}$, and $\overline{\sfCof}$ is the subcategory of all cofibrant objects in $\sfM$ whose morphisms are in $\sfC$.
As an illustration of Theorem \ref{MAIN FOR REP CAT},
we show in Example \ref{fur illusta on model} that
the above two Waldhausen categories are actually coincident.

The paper is organized as follows. In Section \ref{Preliminaries} we introduce necessary notions, examples and facts on Waldhausen categories and Grothendieck opfibrations. Section \ref{proof of mian constru result} is devoted to proving Theorem \ref{total E IS WALD} and revisiting two classical examples due to Waldhausen \cite{Wald83}. The main content of Section \ref{Appl in rep cat} is a proof of Theorem \ref{MAIN FOR REP CAT}.

\section{Preliminaries}\label{Preliminaries}
\noindent In this section, we give necessary preliminaries on Waldhausen categories and
Grothendieck opfibrations.

\subsection{Waldhausen categories}
\label{Gro opfib}
In the original definition of a Waldhausen category \cite{Wald83}, the underlying category is required to be \emph{pointed}. However, some categories we discuss in the paper admit although both initial and terminal objects, they are not necessary isomorphic; see for instance Example \ref{over cate is}. But they still admit classes of cofibrations and weak equivalences satisfying the required axioms. Hence, in this paper, we consider the following slight generalisation of the usual notion of a Waldhausen category, where only an initial object is required.

\begin{dfn}\label{df of wald cat and fun}
A \emph{Waldhausen category} consists of a category $\sfE$ with an initial object ``$0$" and a Waldhausen structure, that is, a class $\sfC$ of morphisms called \emph{cofibrations} (denoted by ``$\rightarrowtail$") and a class $\sfW$ of morphisms called \emph{weak equivalences} (denoted by ``$\overset{\sim}\rightarrow$"), such that both $\sfC$ and $\sfW$ are closed under compositions and the following axioms are satisfied:
\begin{rqm}
\item[(C1)] All isomorphisms in $\sfE$ are cofibrations.
\item[(C2)] For every object $A$ in $\sfE$, the morphism $0 \to A$ is a cofibration.
\item[(C3)] If $A \rightarrowtail B$ is a cofibration, then for any morphism $A \to C$,
            the pushout
                 \[
                   \xymatrix@R=1cm@C=1cm{
                     A
                       \ar@{}[rd]^(0.55)>>{\lrcorner}
                       \ar[d]
                       \ar@{>->}[r]
                   & B
                       \ar[d]                                  \\
                     C
                       \ar[r]
                   & B \sqcup_A C    }
                 \]
            exists in $\sfE$ and the morphism $C \to B \sqcup_A C$ is a cofibration.
\item[(W1)] All isomorphisms in $\sfE$ are weak equivalences.
\item[(W2) ``Gluing Lemma":] Given a commutative diagram
                \[\xymatrix@R=1cm@C=1cm{
                  C
                    \ar[d]^-{\sim}
                & A
                    \ar[l]
                    \ar[d]^-{\sim}
                    \ar@{>->}[r]
                & B
                    \ar[d]^-{\sim}                        \\
                  C'
                & A'
                    \ar[l]
                    \ar@{>->}[r]
                & B'}
                \]
                in $\sfE$, where the horizontal morphisms on the right are cofibrations and all vertical morphisms are weak equivalences, the induced morphism
                \[
                B \sqcup_A C \to B' \sqcup_{A'} C'
                \]
                is also a weak equivalence.
\end{rqm}

A functor between Waldhausen categories is said to be \emph{exact} if it preserves the initial objects, cofibrations, weak equivalences, and the pushout diagrams of axiom (C3).
\end{dfn}

There are many examples of Waldhausen categories. Furthermore, given a Waldhausen category $\sfE$, one can equip Waldhausen structures on a few categories associated to $\sfE$. The following examples are taken from \cite{Wald83}.

\begin{exa}\label{classical exa mor}
Let $\sfE$ be a Waldhausen category.
\begin{rqm}
\item
The morphism category $\Mor(\sfE)$ admits an initial object, and is a Waldhausen category: a morphism $(u, a) : (f: A \rightarrowtail X) \to (g : B \rightarrowtail Y)$ is a cofibration (resp., weak equivalence) in $\Mor(\sfE)$ if both $u$ and $a$ are cofibrations (resp., weak equivalences) in $\sfE$.

\item
Let $\cMor(\sfE)$ be the full subcategory of $\Mor(\sfE)$ whose objects are cofibrations in $\sfE$. For any morphism $(u, a) : (f: A \rightarrowtail X) \to (g : B \rightarrowtail Y)$ in $\cMor(\sfE)$, by considering the following commutative diagram, where the inner square is a pushout that exists by the axiom (C3), we obtain a natural morphism $h$:
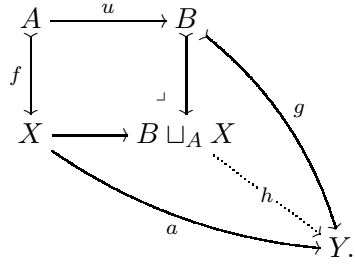
\[
\xymatrix@R=1cm@C=1cm{
   A
     \ar@{>->}[d]_-{f}
     \ar[r]^-{u}
     \ar@{}[rd]^(0.6)>>{\lrcorner}
&  B
     \ar@{>->}[d]_{}
     \ar@{>->}@/^0.8pc/[ddr]^{g}                          \\
   X
     \ar[r]^-{}
     \ar@/_0.8pc/[drr]_{a}
&  B \sqcup_A X
     \ar@{.>}[dr]|-{\, h \,}                              \\
&& Y.}
\]
The category $\cMor(\sfE)$ admits an initial object, and is a Waldhausen category, where $(u, a)$ is a cofibration (resp., weak equivalence) if both $u$ and the induced morphism $h$ are cofibrations (resp., weak equivalences) in $\sfE$.
\end{rqm}
\end{exa}

\begin{exa}\label{over cate is}
Let $\sfE$ be a Waldhausen category with $0$ the initial object and $A$ an object in $\sfE$.
\begin{rqm}
\item
The overcategory $\sfE/A$ admits an initial object $(0\to A)$ as well as a terminal object $\id_A: A\to A$. Thus $\sfE/A$ is not pointed unless $A = 0$. We define a morphism in $\sfE/A$ to be a cofibration (resp., weak equivalence) if it as a morphism in $\sfE$ is a cofibration (resp., weak equivalence). It is routine to check that $\sfE/A$ forms a Waldhausen category with this construction.

\item The undercategory $A/\sfE$ admits an initial object $(\id_A: A\to A)$. However, the construction for overcategories in general does not work for undercategories. Indeed, the unique morphism in $A/\sfE$ from the initial object $(\id_A: A\to A)$ to any object $(f : A \rightarrow X)$ is $f$, which might not be a cofibration in $\sfE$, so in this case the axiom (C2) fails. However, if we consider the full subcategory $\overline{A/\sfE}$ consisting of cofibrations $(A \rightarrowtail X)$ in $\sfE$, then this problem is overcome: one can define a morphism in $\overline{A/\sfE}$ to be a cofibration (resp., weak equivalence) if it as a morphism in $\sfE$ is a cofibration (resp., weak equivalence). With this construction $\overline{A/\sfE}$ forms a Waldhausen category.
\end{rqm}
\end{exa}

\begin{bfhpg}[\bf Weak factorization systems]
Let $l: A \to B$ and $r: C \to D$ be morphisms in a category $\sfE$. We say that $l$ has the {\it left lifting property} with respect to $r$ (or $r$ has the {\it right lifting property} with respect to $l$) if for every pair of morphisms $f: A \to C$ and $g: B \to D$ with
$r \circ f = g \circ l$, there exists a morphism $t: B \to C$ such that the diagram
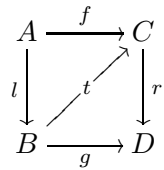
\[
\xymatrix@R=1cm@C=1cm{
  A
\ar[d]_{l}
\ar[r]^{f}
& C
\ar[d]^{r}   \\
  B
\ar[ur]|-{\,t\,}
\ar[r]_{g}
& D
}
\]
commutes.
For a class $\sfC$ of morphisms in $\sfE$, let $\sfC^\Box$ denote the class of morphisms in $\sfE$ having the right lifting property with respect to all morphisms in $\sfC$. The class $^\Box\sfC$ is defined dually. Recall from Bousfield \cite{Bo77} that a pair $(\sfC, \sfF)$ of classes of morphisms in $\sfE$ is called a {\it weak factorization system} if
$\sfC^\Box = \sfF$, ${^\Box\sfF} = \sfC$ and every morphism $\alpha$ in $\sfE$ can be factored as $\alpha = f \circ c$ with $c$ in $\sfC$ and $f$ in $\sfF$.
\end{bfhpg}

\begin{rmk}\label{need pro for wfs}
Let $(\sfC, \sfF)$ be a weak factorization system in $\sfE$.
Then
\begin{prt}
\item
both $\sfC$ and $\sfF$ are closed under compositions, and contain all isomorphisms in $\sfE$;

\item
if $\sfE$ has pushouts along morphisms in $\sfC$, then $\sfC$ is closed under pushouts;

\item
if $\sfE$ has small coproducts, then $\sfC$ is closed under small coproducts.
\end{prt}
For details, see \cite[Propositions D.1.2 and D.1.3]{Jo08}.
\end{rmk}

The following result describes a method to construct a Waldhausen structure from a weak factorization system satisfying certain conditions.

\begin{prp}\label{wfs induce wald}
Let $\sfE$ be a category with finite coproducts and an initial object 0, and let $(\sfC, \sfW)$ be a pair of classes of morphisms in $\sfE$ satisfying the following conditions:
\begin{rqm}
\item $(\sfC, \sfC^\Box)$ forms a weak factorization system in $\sfE$ such that
      $\sfC^\Box \subseteq \sfW$;

\item $\sfW$ contains all isomorphisms in $\sfE$ and satisfies the 2-out-of-3 property, that is, if two of the three morphisms $f$, $g$ and $gf$ are in $\sfW$, then so is the third one;

\item $\sfC \cap \sfW$ is closed under pushouts;

\item $0 \to X \in \sfC$ for all objects $X$ in $\sfE$.
\end{rqm}
Then $(\sfC, \sfW)$ forms a Waldhausen structure on $\sfE$.
\end{prp}

\begin{prf*}
By Remark \ref{need pro for wfs}, $\sfC$ contains all isomorphisms in $\sfE$ and is closed under pushouts. It remains to check the Gluing Lemma. This could be done by showing that $\sfE$ is a category of cofibrant objects in the sense of Kamps and Porter \pagecite[80]{K-P97}, as any category of cofibrant objects satisfies the Gluing Lemma by \thmcite[II 2.27]{K-P97}. However, according to the axioms for $\sfE$ being a category of cofibrant objects, we only need to show that any object $X \in \sfE$ admits a cylinder object. Indeed, since $(\sfC, \sfC^\Box)$ forms a weak factorization system in $\sfE$, we can factor the fold map $\id_X \amalg \id_X : X \amalg X \to X$ as
\[
X \amalg X \overset{e_0 \amalg e_1} \longrightarrow
X \times I \overset{\delta}         \longrightarrow
X
\]
for some object $X \times I$ in $\sfE$ with $e_0 \amalg e_1 \in \sfC$ and $\delta \in \sfC^\Box$. But $\sfC^\Box \subseteq \sfW$, so $\delta \in \sfW$, and hence $(X \times I, e_0, e_1, \delta)$ is the desired cylinder object for $X$.
\end{prf*}

We illustrate the application of this proposition in the following examples.

\begin{exa}\label{exa model}
Let $\sfM$ be a \emph{Quillen model category}; that is, it has finite limits and colimits, and there is a triple $(\sfC, \sfW, \sfF)$ of classes of morphisms in $\sfM$ such that both $(\sfC, \sfW \cap \sfF)$ and $(\sfC \cap \sfW, \sfF)$ are weak factorization systems, and $\sfW$ satisfies the 2-out-of-3 property. 
Recall that an object $M$ in $\sfM$ is called \emph{cofibrant} if the morphism $0 \to M$ belongs to $\sfC$, where $0$ is the initial object of $\sfM$. By Proposition \ref{wfs induce wald}, the full subcategory $\sfCof$ of $\sfM$ consisting of all cofibrant objects forms a Waldhausen category with $(\Mor(\sfCof) \cap \sfC, \, \Mor(\sfCof) \cap \sfW)$ its Waldhausen structure.
\end{exa}

The following example, which is essentially due to Sarazola \cite{Sarazola20}, shows that a complete cotorion pair (see \cite{rha} for a definition) in an abelian category $\sfA$ might induce a Waldhausen structure on $\sfA$. For a full subcategory $\sfX$ of $\sfA$, set
\begin{align*}
\Mon(\sfX) & = \{\alpha\ |\
                 \alpha\ \text {is a monomorphism}\ \text{with}\
                 \Coker(\alpha) \in \sfX\}, \ \textrm{and} \\
\Epi(\sfX) & = \{\alpha\ |\
                 \alpha\ \text {is an epimorphism}\ \text{with}\
                 \Ker(\alpha) \in \sfX\}.
\end{align*}

\begin{exa}\label{complete cotor pair}
Let $\sfA$ be an abelian category with enough injectives, and denote by $\sfInj$ (resp., $\sfInj^{<\infty}$) the full subcategory of $\sfA$ consisting of all injectives (resp., objects with finite injective dimension). We use Proposition \ref{wfs induce wald} to show that the pair $(\Mon(\sfA), \sfW)$ gives a Waldhausen structure on $\sfA$, where $\sfW$ is the class of morphisms $f$ that can be factored as $f = g \circ h$ with $h \in \Mon(\sfInj^{<\infty})$ and $g \in \Epi(\sfInj)$.

We check that $(\Mon(\sfA), \sfW)$ satisfies conditions (1)-(4) in Proposition \ref{wfs induce wald}. For (1), since $(\sfA, \sfInj)$ is a complete cotorsion pair in $\sfA$, it follows from \thmcite[2.4]{Ho02} that $(\Mon(\sfA), \Epi(\sfInj))$ forms a weak factorization system in $\sfA$. For (2), it is obvious that $\sfW$
contains all isomorphisms in $\sfA$. Since $\sfInj^{<\infty}$ satisfies the 2-out-of-3 property for short exact sequences in $\sfA$, we conclude by \prpcite[5.6]{Sarazola20} that $\sfW$ satisfies the 2-out-of-3 property for composable morphisms in $\sfA$. For (3), it is easy to see that $\Mon(\sfA) \cap \sfW = \Mon(\sfInj^{<\infty})$. Thus, $\Mon(\sfA) \cap \sfW$ is closed under pushouts as so is $\Mon(\sfInj^{<\infty})$. Finally, $(\Mon(\sfA), \sfW)$ clearly satisfies (4).
\end{exa}

\subsection{Grothendieck opfibrations}\label{Gro opfib}
Let $p : \calT \to \calB$ be a functor. We say that an object $X$ in $\calT$ is \emph{above} an object $A$ in $\calB$ if $p(X) = A$ and, similarly, that a morphism $f : X \to Y$ in $\calT$ is \emph{above} a morphism $u: A \to B$ in $\calB$ if $p(f) = u$. Given an object $A$ in $\calB$, the \emph{fiber} of $A$ with respect to $p$ is the full subcategory of $\calT$ whose objects are those $X$ in $\calT$ such that $p(X)=A$ and morphisms are those $f$ in $\calT$ such that $p(f)=\id_A$. We denote the fiber of $A$ with respect to $p$ by $\calT_A$.

Recall that a morphism $f : X \to Y$ in $\calT$ that is above a morphism $u : A \to B$ in $\calB$ is called \emph{cocartesian} with respect to $p$ if for every pair of morphisms $v : B \to C$ in $\calB$ and $g : X \to Z$ in $\calT$ that is above $v \circ u : A \to C$, there exists a unique morphism $h : Y \to Z$ such that it is above $v$ and $h \circ f = g$. Diagrammatically:
\[
\xymatrix@R=0.4cm@C=0.5cm{
&&&&& Z
\ar@{.}[dd]
\\
      X
\ar[rrr]_{f}
\ar@/^0.85pc/[rrrrru]^{g}
\ar@{.}[dd]
&&&   Y
\ar@.[rru]_{h}
\ar@{.}[dd]                          \\
&&& \ar@/^0.2pc/[rr] && C            \\
      A
\ar[rrr]_{u}
\ar@/^0.3pc/@{}[rrru]^(0.85){v \circ u}
&&&   B
\ar[rru]_{v}   }
\]
It is easy to check that the class of all cocartesian morphisms is closed under compositions.

As an immediate consequence of the universal property of cocartesian morphisms, we have:

\begin{lem}\label{unique for cocar mor}
Let $f : X \to Y$  be a cocartesian morphism and $g, g' : Y \to Z$ morphisms with
$g \circ f = g' \circ f$. If $p(g) = p(g')$, then one has $g = g'$.
\end{lem}

\begin{dfn}\label{df of Grothendieck opfibra}
A \emph{cloven Grothendieck opfibration} is a functor $p : \calT \to \calB$ together with for any object $X$ in $\calT$ that is above an object $A$ in $\calB$ and any morphism $u : A \to B$ in $\calB$,
\begin{itemize}
\item
an object $u_!(X)$ in $\calT$ and a cocartesian morphism $\lambda^p_{u,X} : X \to u_!(X)$ in $\calT$ that is above $u$; diagrammatically:
\[
\xymatrix@R=0.5cm@C=1cm{
   X
\ar[r]^-{\lambda^p_{u, X}}
\ar@{.}[d]
&  u_!(X)
\ar@{.}[d]           \\
   A
\ar[r]^-{u}
&  B }
\]
\end{itemize}
When the functor $p$ is clear from the context, we write $\lambda_{u,X}$ for $\lambda^p_{u, X}$, and call it the \emph{cocartesian lifting} of $X$ along $u$.
\end{dfn}

{\bf In this paper}, we always consider cloven Grothendieck opfibrations, and call them simply Grothendieck opfibrations.
Indeed, if $\calT$ and $\calB$ are small relatively to a universe $\mathbb{U}$ in which the axiom of choice is assumed,
then a cloven Grothendieck opfibration is exactly the same as the original notion of a Grothendieck opfibration; see for instance \cite{MoVa20}.



The following result will be used frequently in the paper.

\begin{lem}\label{key factorization}
Let $p: \calT \to \calB$ be a Grothendieck opfibration, and $f: X \to Y$ a morphism in $\calT$ that is above a morphism $u: A\to B$ in $\calB$. Then $f$ can be factored uniquely as the cocartesian morphism $\lambda_{u, X}: X \to u_!(X)$ followed by a morphism $f_\triangleright: u_!(X)\to Y$ in the fiber $\calT_B$. If furthermore $f$ is an isomorphism, then $f_\triangleright$ is an isomorphism in $\calT_B$.
\end{lem}

\begin{prf*}
Consider the following diagram of solid arrows:
\begin{equation*} \label{factorization 1}
\tag{\ref{key factorization}.1}
\begin{gathered}
\xymatrix@R=0.4cm@C=0.5cm{
&&&&&
Y
\ar@{.}[dd]
\\
X
\ar[rrr]_{\lambda_{u, X}}
\ar@/^0.85pc/[rrrrru]^{f}
\ar@{.}[dd]
&&&
u_!(X)
\ar@{.>}[rru]_{f_\triangleright}
\ar@{.}[dd]
\\
&&& \ar@/^0.2pc/[rr] &&
B
\\
A
\ar[rrr]_{u}
\ar@/^0.3pc/@{}[rrru]^(0.85){u}
&&&
B
\ar@{=}[rru]}
\end{gathered}
\end{equation*}
By the cocartesian universal property of $\lambda_{u, X}$, there is a unique morphism $f_\triangleright$, which is clearly in the fiber $\calT_B$, such that $f_\triangleright\circ\lambda_{u,X}=f$.

If $f$ is an isomorphism, let $g: Y \to X$ be the inverse of $f$. We would like to find the inverse of $f_\triangleright$. To this end, set $v = p(g): B \to A$. By the above proof, $g$ can be factored uniquely as the cocartesian morphism $\lambda_{v, Y}: Y \to v_!(Y)$ followed by a morphism $g_\triangleright: v_!(Y)\to X$ in the fiber $\calT_A$. Note that
\[ (\lambda_{u,X} \circ g_\triangleright  \circ \lambda_{v,Y}
                    \circ f_\triangleright) \circ \lambda_{u,X}
= \lambda_{u,X}  \circ g \circ f
= \lambda_{u,X}
= \id_{u_!(X)}  \circ \lambda_{u,X};\]
diagrammatically:
\[
\xymatrix@R=0.75cm@C=1cm{\\
   X \ar[rr]^-{f} \ar[dr]_-{\lambda_{u,X}} \ar@{=}@/^3pc/[rrrr]
&& Y \ar[rr]^-{g} \ar[dr]_-{\lambda_{v,Y}} \ar@{=}@/^3pc/[rrrr]
&& X \ar[rr]^-{f} \ar[dr]_-{\lambda_{u,X}}
&& Y. \\
&  u_!(X) \ar[ur]_-{f_\triangleright}
&& v_!(Y) \ar[ur]_-{g_\triangleright}
&& u_!(X) \ar[ur]_-{f_\triangleright}
}
\]
Moreover, one has
\[p(\lambda_{u,X} \circ g_\triangleright \circ \lambda_{v,Y} \circ  f_\triangleright)
= u \circ v \circ \id_{B} = p(fg) = p(\id_Y) = \id_{B}= p(\id_{u_!(X)}).\]
Since $\lambda_{u,X}$ is cocartesian, it follows from Lemma \ref{unique for cocar mor} that
\[
(\lambda_{u,X} \circ g_\triangleright \circ\lambda_{v,Y}) \circ  f_\triangleright
= \id_{u_!(X)}.
\]
It is also clear that
$f_\triangleright \circ (\lambda_{u, X} \circ g_\triangleright \circ\lambda_{v, Y})
=f\circ g= \id_Y.$
Thus, $\lambda_{u, X} \circ g_\triangleright \circ\lambda_{v,Y}$ is the inverse of $f_\triangleright$.
\end{prf*}

\begin{rmk}\label{the mor u!(k)}
By this lemma, the $u_!$ appearing in Definition \ref{df of Grothendieck opfibra} is a functor from $\calT_A$ to $\calT_B$, called the \emph{reindexing functor} of $u$. To  specify its action on morphisms, let $k: X\to X'$ be a morphism in $\calT_A$, and consider the morphism $\lambda_{u,X'}\circ k: X \to u_!(X')$,
where $\lambda_{u,X'}: X' \to u_!(X')$ is the cocartesian lifting of $X'$ along $u$. Then by Lemma \ref{key factorization}, one has $\lambda_{u,X'}\circ k=(\lambda_{u,X'}\circ k)_\triangleright\circ\lambda_{u,X}$, so the following diagram commutes:
\begin{equation*} \label{factorization 2}
\tag{\ref{the mor u!(k)}.1}
\begin{gathered}
\xymatrix@R=1cm@C=1cm{
  X
   \ar[d]_-{k} \ar[r]^-{\lambda_{u,X}}
& u_!(X)
   \ar[d]^-{(\lambda_{u,X'}\circ k)_\triangleright} \\
  X'
   \ar[r]_-{\lambda_{u,X'}}
& u_!(X')                           }
\end{gathered}
\end{equation*}
Thus we can set $u_!(k)=(\lambda_{u,X'} \circ k)_\triangleright$. It is easy to check that this construction is functorial.
\end{rmk}

\begin{rmk}\label{assume pseudo to be functor}
The rule assigning an object $A$ in $\calB$ to the fiber $\calT_A$ and a morphism $u: A\to B$ in $\calB$ to the reindexing functor
$u_!: \calT_A \to \calT_B$ gives a pseudo functor $\scrF: \calB \to \Cat$. For two morphisms $u: A \to B$ and $v: B\to C$ in $\calB$, there is a natural isomorphism
$\phi: (v \circ u)_! \to v_! \circ u_!$. Thus, for any object $X$ in $\calT_A$, there is an isomorphism $\phi_X: (v \circ u)_!(X) \to v_!(u_!(X))$ in $\calT_C$ such that $\lambda_{v,u_!(X)} \circ \lambda_{u,X}=  \phi_X \circ \lambda_{v \circ u,X}$. To simplify subsequent proofs, {\bf from now on,} we identify $(v \circ u)_!$ with $v_! \circ u_!$ via the isomorphism $\phi_X$. Accordingly, the above identity can be simplified to
\begin{equation*} \label{useful equ}
\tag{\ref{assume pseudo to be functor}.1}
\lambda_{v,u_!(X)}\circ\lambda_{u,X}=\lambda_{v \circ u,X}.
\end{equation*}
\end{rmk}

The following result shows, under certain mild conditions, that the total category $\calT$ and the fiber $\calT_{0}$ of the initial object $0$ in the basis category $\calB$ share the common initial object, and that $\calT$ is pointed if so are $\calB$ and $\calT_{0}$.

\begin{lem}\label{intial object}
Let $p: \calT \to \calB$ be a Grothendieck opfibration.
Suppose that the basis category $\calB$ admits an initial object $0$ and the fiber $\calT_{0}$ admits an initial object $\overline{0}$.
If the reindexing functor $u_!$ preserves the initial objects for each morphism $u$ in $\calB$,
then $\overline{0}$ is the initial object of the total category $\calT$. Furthermore, if both $\calB$ and $\calT_{0}$ are pointed categories, then so is $\calT$.
\end{lem}

\begin{prf*}
Let $X$ be an object in $\calT$.
Consider the diagram
\[
\xymatrix@R=0.5cm@C=0.75cm{
\overline{0}
\ar[rr]^-{\lambda_{\ast, \overline{0}}}
\ar@{.}[d]
& & \ast_!(\overline{0})
\ar@{.}[d]
\ar[rr]^-{\exists \, \mid \,\, \overline{\ast}}
& & X
\ar@{.}[d]  \\
0
\ar[rr]^-{\exists \, \mid \,\, \ast}
& & p(X)
\ar@{=}[rr]
& & p(X)}
\]
in which $\lambda_{\ast, \overline{0}}$
is the cocartesian lifting of $\overline{0}$ along $\ast$, and
$\ast_!(\overline{0})$ is the initial object in $\calT_{p(X)}$ as
$\ast_!: \calT_{0}\to \calT_{p(X)}$ preserves the initial objects by assumption.
Then we obtain a morphism
$\overline{\ast} \circ \lambda_{\ast, \overline{0}}: \overline{0} \to X$; we show next that it is the unique morphism from $\overline{0}$ to $X$. Let $f: \overline{0} \to X$ be another morphism in $\calT$. Then $f$ is clearly above $\ast$, and so by
Lemma \ref{key factorization}, one has $f=f_\triangleright\circ\lambda_{\ast, \overline{0}}$ with $f_\triangleright: \ast_!(\overline{0})\to X \in \calT_{p(X)}$. We mention that $\ast_!(\overline{0})$ is the initial object in $\calT_{p(X)}$. It follows that
$f_\triangleright = \overline{\ast}$, and so
$f = f_\triangleright \circ \lambda_{\ast, \overline{0}}
= \overline{\ast} \circ \lambda_{\ast, \overline{0}}$,
as desired. Thus, $\overline{0}$ is the initial object in $\calT$.

If $0$ is a terminal object in $\calB$ and $\overline{0}$ is a terminal object in $\calT_{0}$, then using a similar proof as above, one gets that $\overline{0}$ is a terminal object in $\calT$ as well.
\end{prf*}

\section{Proof of Theorem \ref{total E IS WALD} and two classical examples}
\label{proof of mian constru result}
\noindent This section is devoted to giving a proof of Theorem \ref{total E IS WALD} and revisiting two classical examples due to Waldhausen \cite{Wald83} as illustrations of Theorem \ref{total E IS WALD}, where some details will be applied further in the proof of Theorem \ref{MAIN FOR REP CAT}.

\subsection{Proof of Theorem \ref{total E IS WALD}}
Recall that a Grothendieck opfibration $p: \calT \to \calB$ is called
a \emph{Waldhausen opfibration} if the following conditions are satisfied:
\begin{itemize}
\item for each object $A$ in $\calB$,
      the fiber $\calT_A$ is a Waldhausen category;
\item for each morphism $u: A\to B$ in $\calB$,
      the reindexing functor $u_!: \calT_A\to \calT_B$ is exact, that is, $u_!$ preserves the initial objects, cofibrations, weak equivalences, and the pushout diagrams of axiom (C3).
\end{itemize}
{\bf Throughout this subsection}, we always assume that $p: \calT \to \calB$ is a Waldhausen opfibration with the basis category $\calB$ a Waldhausen category. Recall that a \emph{total cofibration} (resp., \emph{total weak equivalence}) is a morphism $f: X \to Y$ in $\calT$ that is above a cofibration (resp., \emph{weak equivalence}) $u: A \to B$ in $\calB$ such that $f_\triangleright$ is a cofibration (resp., \emph{weak equivalence}) in the fiber $\calT_B$. We aim to prove that the classes of total cofibrations and total weak equivalences
endow the category $\calT$ with a Waldhausen structure.

\begin{lem}\label{both are subcat}
Both the class of total cofibrations and the class of total weak equivalences are closed under compositions.
\end{lem}

\begin{prf*}
We only prove the conclusion for total cofibrations. A similar proof works for total weak equivalences.

Let $X \overset{f} \longrightarrow Y \overset{g} \longrightarrow Z$ be a pair of total cofibrations in $\calT$. Consider the picture
\[
\xymatrix@R=1cm@C=1cm{
& &  Z
\\
& Y \ar[ru]^-{g} \ar[r]_-{\lambda_{v,Y}}
& v_!(Y) \ar[u]_-{g_\triangleright}
\\
X \ar[ru]^-{f} \ar[r]_-{\lambda_{u,X}}
& u_!(X) \ar[r]_-{\lambda_{v,u_!(X)}} \ar[u]_-{f_\triangleright}
& v_!(u_!(X))\ar[u]_-{v_!(f_\triangleright)}
}
\xymatrix@R=0.75cm@C=1cm{ \\ \quad \overset{p}\longmapsto \quad}
\xymatrix@R=1cm@C=1.25cm{
& &  C
\\
& B \ar[ru]^-{v} \ar[r]_-{v}
& C \ar@{=}[u]_-{}
\\
A \ar[ru]^-{u} \ar[r]_-{u}
& B \ar[r]_-{v} \ar@{=}[u]
& C\ar@{=}[u]
}
\]
where the left commutative diagram in $\calT$ is above the right commutative diagram in $\calB$, and the inner square in the left diagram is from (\ref{factorization 2}).

To show that $g \circ f$ is also a total cofibration, we need to prove that $v \circ u$ is a cofibration in $\calB$ and $(g \circ f)_\triangleright$ is a cofibration in $\calT_C$. Since $f$ and $g$ are total cofibrations, it follows that $u$ and $v$ are cofibrations in $\calB$, $f_\triangleright$ is a cofibration in $\calT_B$,
and $g_\triangleright$ is a cofibration in $\calT_C$. Thus $v \circ u$ is a cofibration in $\calB$. Furthermore, we have $(g \circ f)_\triangleright = g_\triangleright \circ v_!(f_\triangleright)$ by equality (\ref{useful equ}). Since $v_!$ is exact, it preserves cofibrations. Hence, $v_!(f_\triangleright)$ is a cofibration in $\calT_C$. Consequently, $(g \circ f)_\triangleright$ is a cofibration in $\calT_C$.
\end{prf*}

By Lemma \ref{intial object}, the total category $\calT$ admits an initial object, and it follows from Lemma \ref{key factorization} that $\calT$ satisfies the axioms (C1) and (W1). In the rest of this subsection we prove that $\calT$ satisfies the axioms (C2), (C3) and (W2).

\begin{bfhpg}[\bf For the axiom (C2)]\label{prove C2}
By Lemma \ref{intial object}, the initial object in $\calT$ is $\overline{0}$, the initial object in the fiber $\calT_{0}$, where $0$ is the initial object in $\calB$. Moreover,  for any object $X$ in $\calT$, the unique morphism $f: \overline{0}\to X$ that is above the morphism $\ast: 0 \to p(X)$ in $\calB$ factors uniquely as
\[
\xymatrix@R=0.5cm@C=0.5cm{
    \overline{0}
      \ar[rr]^-{\lambda_{\ast, \overline{0}}}
& & \ast_!(\overline{0})
      \ar[rr]^-{\exists \, \mid \,\, \overline{\ast}}
& & X,}
\]
where $\ast_!(\overline{0})$ is the initial object in $\calT_{p(X)}$ and $\overline{\ast}=f_\triangleright$. Since both $\calB$ and $\calT_{p(X)}$ are Waldhausen categories, it follows that $\ast$ is a cofibration in $\calB$ and $f_\triangleright$ is a cofibration in $\calT_{p(X)}$. Thus, $f$ is a total cofibration in $\calT$. \qed
\end{bfhpg}


\begin{bfhpg}[\bf For the axiom (C3)]\label{prove C3}
Let $\xymatrix{Z & X \ar[l]_-{g} \ar@{>->}[r]^-{f} & Y}$
be a pair of morphisms in $\calT$ with $f$ a total cofibration.
It is sufficient to construct the pushout of $f$ along $g$ and
show that it preserves total cofibrations. Suppose that the above pair of morphisms is above the one $\xymatrix{C& A \ar[l]_-{v} \ar@{>->}[r]^-{u} & B}$
in $\calB$. Then $u$ is a cofibration as $f$ is a total cofibration. Since $\calB$ satisfies the axiom (C3),
there exists a pushout
\begin{equation*} \label{pushout 1}
\tag{\ref{prove C3}.1}
\begin{gathered}
\xymatrix@R=1cm@C=1cm{
  A
  \ar[d]_-{v}
  \ar@{}[rd]^(0.66)>>{\lrcorner}
  \ar@{>->}[r]^-{u}
& B
  \ar[d]^-{\overline{v}}   \\
  C
  \ar@{>->}[r]^-{\overline{u}}
& B\sqcup_AC                        }
\end{gathered}
\end{equation*}
with $\overline{u}$ a cofibration in $\calB$.
We mention that $f_\triangleright: u_!(X) \to Y$ is a cofibration in $\calT_B$.
It follows that $\overline{v}_!(f_\triangleright)$ is a cofibration in $\calT_{B\sqcup_AC}$,
as $\overline{v}_!$ is exact by the assumption. Since $\calT_{B\sqcup_AC}$ satisfies the axiom (C3) as well,
there exists the pushout
\begin{equation*} \label{pushout 2}
\tag{\ref{prove C3}.2}
\begin{gathered}
\xymatrix@R=1cm@C=1cm{
  \overline{u}_!(v_!(X)) = \overline{v}_!(u_!(X))
    \ar[d]_-{\overline{u}_!(g_\triangleright)}
    \ar@{}[rd]^(0.66)>>{\lrcorner}
    \ar@{>->}[r]^-{\overline{v}_!(f_\triangleright)}
&  \overline{v}_!(Y) \ar[d]^-{a}   \\
  \overline{u}_!(Z)
    \ar@{>->}[r]^-{b}
& W                          }
\end{gathered}
\end{equation*}
with $b$ a cofibration in $\calT_{B\sqcup_AC}$.
This yields the fourth equality in the following computation:
\begin{align*}
   a \circ \lambda_{\overline{v}, Y}
     \circ f
&= a \circ \lambda_{\overline{v}, Y}
     \circ f_\triangleright
     \circ \lambda_{u, X}
& \text{by (\ref{factorization 1})}                                               \\
&= a \circ \overline{v}_!(f_\triangleright)
     \circ \lambda_{\overline{v},u_!(X)}
     \circ \lambda_{u, X}
& \text{by (\ref{factorization 2}) for $f_\triangleright$}                                               \\
&= a \circ \overline{v}_!(f_\triangleright)
     \circ \lambda_{\overline{u},v_!(X)}
     \circ \lambda_{v, X}
& \text{by (\ref{useful equ}) as }
        \overline{u} \circ v = \overline{v} \circ u                               \\
&= b \circ \overline{u}_!(g_\triangleright)
     \circ \lambda_{\overline{u},v_!(X)}
     \circ \lambda_{v, X}                                                         \\
&= b \circ \lambda_{\overline{u}, Z}
     \circ g_\triangleright
     \circ \lambda_{v, X}
& \text{by (\ref{factorization 2}) for $g_\triangleright$}                                               \\
&= b \circ \lambda_{\overline{u}, Z}
     \circ g
& \text{by (\ref{factorization 1}) }
\end{align*}
Thus the square
\begin{equation*} \label{pushout 3}
\tag{\ref{prove C3}.3}
\begin{gathered}
\xymatrix@R=1cm@C=1.5cm{
   X \ar@{>->}[r]^-{f} \ar[d]_-{g}
&  Y \ar[d]^-{a \circ \lambda_{\overline{v}, Y}} \\
   Z \ar[r]^-{b \circ \lambda_{\overline{u}, Z}}
&  W }
\end{gathered}
\end{equation*}
in $\calT$ commutes. Since $b \circ \lambda_{\overline{u}, Z}$ is above $\overline{u}$ and $(b \circ \lambda_{\overline{u}, Z})_\triangleright = b$, according to what we established above, $b \circ \lambda_{\overline{u}, Z}$ is a total cofibration.

It remains to prove that (\ref{pushout 3}) is a pushout. As depicted below, the commutative diagram of solid arrows in $\calT$ is above the commutative diagram of solid arrows in $\calB$:
\[
\xymatrix@R=1cm@C=1.5cm{
   X
   \ar@{>->}[r]^-{f}
   \ar[d]_-{g}
&  Y
   \ar[d]_-{a \circ \lambda_{\overline{v}, Y}}
   \ar@/^/[ddr]^{c}                                       \\
   Z
   \ar[r]^-{b \circ \lambda_{\overline{u}, Z}}
   \ar@/_/[drr]_{d}
&  W
   \ar@{.>}[dr]|-{e}                                      \\
&& V }
\xymatrix@R=1cm@C=1cm{ \\ \quad \overset{p}\longmapsto \quad}
\xymatrix@R=1cm@C=1cm{
   A
   \ar[d]_-{v}
   \ar@{>->}[r]^-{u}
   \ar@{}[rd]^(0.6)>>{\lrcorner}
&  B
   \ar[d]_{\overline{v}}
   \ar@/^/[ddr]^{p(c)}                                     \\
   C
   \ar@{>->}[r]^-{\overline{u}}
   \ar@/_/[drr]_{p(d)}
&  B\sqcup_AC
   \ar@{.>}[dr]|-{t}                                       \\
&& p(V)              }
\]
We need to find a unique morphism $e$ such that $e \circ (b \circ \lambda_{\overline{u}, Z}) = d$ and $e \circ (a \circ \lambda_{\overline{v}, Y}) = c$.

Since (\ref{pushout 1}) is a pushout, there exists a unique morphism $t$ such that $t \circ \overline{u} = p(d)$ and $t \circ \overline{v} = p(c)$. Since (\ref{pushout 2}) is a pushout in $\calT_{B\sqcup_AC}$ and $t_!$ preserves the pushout of axiom (C3), the square in the following diagram is a pushout in $\calT_{p(V)}$:
\begin{equation*} \label{pushout 4}
\tag{\ref{prove C3}.4}
\begin{gathered}
\xymatrix@R=1cm@C=1.5cm{
    t_!(\overline{u}_!(v_!(X))) = t_!(\overline{v}_!(u_!(X)))
      \ar[d]_-{t_!(\overline{u}_!(g_\triangleright))}
      \ar@{}[rd]^(0.6)>>{\lrcorner}
      \ar@{>->}[r]^-{t_!(\overline{v}_!(f_\triangleright))}
&   t_!(\overline{v}_!(Y)) = p(c)_!(Y)
      \ar[d]_{t_!(a)}
      \ar@/^/[ddr]^{c_\triangleright}                              \\
    p(d)_!(Z) = t_!(\overline{u}_!(Z))
      \ar@{>->}[r]^-{t_!(b)}
      \ar@/_/[drr]_{d_\triangleright}
&   t_!(W)
      \ar@{.>}[dr]|-{h}                                            \\
&&  V  }
\end{gathered}
\end{equation*}
However, we have
\begin{align*}
&  \quad\,\,
    c_\triangleright \circ t_!(\overline{v}_!(f_\triangleright))
                     \circ \lambda_{p(c), u_!(X)}
                     \circ \lambda_{u, X} \\
& = c_\triangleright \circ p(c)_!(f_\triangleright)
                     \circ \lambda_{p(c), u_!(X)}
                     \circ \lambda_{u, X}
& \text{ as } t \circ \overline{v} = p(c) \\
& = c_\triangleright \circ \lambda_{p(c), Y}
                     \circ f_\triangleright
                     \circ \lambda_{u, X}
& \text{by (\ref{factorization 2}) for $f_\triangleright$} \\
& = c \circ f
& \text{by (\ref{factorization 1})}
\end{align*}
and
\begin{align*}
&  \quad\,\,
    d_\triangleright \circ t_!(\overline{u}_!(g_\triangleright))
                     \circ \lambda_{p(c), u_!(X)}
                     \circ \lambda_{u, X}                                           \\
& = d_\triangleright \circ p(d)_!(g_\triangleright)
                     \circ \lambda_{p(c), u_!(X)}
                     \circ \lambda_{u, X}
& \text{ as } t \circ \overline{u} = p(d) \\
& = d_\triangleright \circ p(d)_!(g_\triangleright)
                     \circ \lambda_{p(d), v_!(X)}
                     \circ \lambda_{v, X}
& \text{by (\ref{useful equ}) as } p(c) \circ u = p(d) \circ v \\
& = d_\triangleright \circ \lambda_{p(d), Z}
                     \circ g_\triangleright
                     \circ \lambda_{v, X}
& \text{by (\ref{factorization 2}) for $g_\triangleright$}                                                  \\
& = d \circ g
& \text{by (\ref{factorization 1})}
\end{align*}
Consequently, we have
\[(c_\triangleright \circ t_!(\overline{v}_!(f_\triangleright)))
                    \circ (\lambda_{p(c), u_!(X)} \circ \lambda_{u, X})
= (d_\triangleright \circ t_!(\overline{u}_!(g_\triangleright)))
                    \circ (\lambda_{p(c), u_!(X)} \circ \lambda_{u, X}).
\]
Since $\lambda_{p(c), u_!(X)} \circ \lambda_{u, X}$ is cocartesian and
\[
p(c_\triangleright \circ t_!(\overline{v}_!(f_\triangleright)))
= p(c_\triangleright \circ p(c)_!(f_\triangleright))
= \id_{p(V)} \circ \id_{p(V)}
= p(d_\triangleright \circ p(d)_!(g_\triangleright))
= p(d_\triangleright \circ t_!(\overline{u}_!(g_\triangleright))),
\]
it follows from Lemma \ref{unique for cocar mor} that $c_\triangleright \circ t_!(\overline{v}_!(f_\triangleright)) = d_\triangleright \circ t_!(\overline{u}_!(g_\triangleright))$, so the outermost square in (\ref{pushout 4}) is commutative, and hence there exists a unique morphism $h$ such that
$h \circ t_!(b) = d_\triangleright$ and $h \circ t_!(a) = c_\triangleright$.

Set $e = h \circ \lambda_{t, W}$. We check that $e \circ (b \circ \lambda_{\overline{u}, Z}) = d$ and $e \circ (a \circ \lambda_{\overline{v}, Y}) = c$. Indeed, we have
\begin{align*}
&  \quad\,\,
    e \circ (b \circ \lambda_{\overline{u}, Z})                    \\
& = h \circ \lambda_{t, W}
      \circ b
      \circ \lambda_{\overline{u}, Z}                              \\
& = h \circ t_!(b)
      \circ \lambda_{t, \overline{u}_!(Z)}
      \circ \lambda_{\overline{u}, Z}
& \text{by (\ref{factorization 2}) for $b$}                                \\
& = d_\triangleright \circ \lambda_{t, \overline{u}_!(Z)}
                     \circ \lambda_{\overline{u}, Z}\\
& = d_\triangleright \circ \lambda_{p(d), Z}
& \text{by (\ref{useful equ}) as }
        t \circ \overline{u} = p(d)                                \\
& = d
& \text{by (\ref{factorization 1})}
\end{align*}
Similarly, we have $e \circ (a \circ \lambda_{\overline{v}, Y}) = c$. This completes the proof. \qed
\end{bfhpg}

\begin{bfhpg}[\bf For the axiom (W2)]\label{prove W2}
Let
\[
\xymatrix@R=1cm@C=1cm{
  Z  \ar[d]^-{i}_-{\sim}
& X  \ar[l]_-{g}
     \ar[d]^-{j}_-{\sim}
     \ar@{>->}[r]^-{f}
& Y  \ar[d]^-{k}_-{\sim}            \\
  Z'
& X' \ar[l]_-{g'}
     \ar@{>->}[r]^-{f'}
& Y' }
\]
be a commutative diagram in $\calT$ such that $f$ and $f'$ are total cofibrations and $i$, $j$ and $k$ are total weak equivalences. By \ref{prove C3}, we obtain pushouts in $\calT$ appearing as the top and bottom squares in the following commutative diagram in $\calT$ on the left side (we still adopt the notation used in \ref{prove C3}, although there will be some adjustments that the readers are well aware of).
\begin{equation*} \label{WEAK 1}
\tag{\ref{prove W2}.1}
\begin{gathered}
\xymatrix@R=1cm@C=1cm{
& X
     \ar[dl]_-{g}
     \ar@{>->}[rr]^(0.75){f}
     \ar'[d]^-{j}_-{\sim}[dd]
&  & Y
     \ar[dd]^(0.33){k}_(0.33){\sim}
     \ar[dl]_-{a \circ \lambda_{\overline{v}, Y}}                              \\
Z    \ar@{>->}[rr]^(0.72){b \circ \lambda_{\overline{u}, Z}}
     \ar[dd]_(0.33){i}^(0.33){\sim}
&  & W
     \ar[dd]^(0.33){e}                                  \\
& X' \ar[dl]_-{g'}
     \ar@{>->}'[r][rr]^-{f'}
&  & Y'
     \ar[dl]^-{a' \circ \lambda_{\overline{v'}, Y'}}                          \\
Z'
     \ar@{>->}[rr]_(0.66){b' \circ \lambda_{\overline{u'}, Z'}}
&  & W'       }
\xymatrix@R=1cm@C=0.75cm{\\ \quad \overset{p}\longmapsto \quad}
\xymatrix@R=1cm@C=0.8cm{
& A
     \ar[dl]_-{v}
     \ar@{>->}[rr]^(0.75){u}
     \ar'[d]^-{p(j)}_-{\sim}[dd]
&  & B
     \ar[dd]^(0.33){p(k)}_(0.33){\sim}
     \ar[dl]^-{\overline{v}}                              \\
C
     \ar@{>->}[rr]^(0.66){\overline{u}}
     \ar[dd]_(0.33){p(i)}^(0.33){\sim}
&  & B\sqcup_AC
     \ar[dd]^(0.33){t}                                  \\
& A'
     \ar[dl]_-{v'}
     \ar@{>->}'[r][rr]^-{u'}
&  & B'
     \ar[dl]^-{\overline{v'}}                          \\
C'
     \ar@{>->}[rr]_(0.56){\overline{u'}}
&  & B'\sqcup_{A'}C'       }
\end{gathered}
\end{equation*}
We would like to show that the induced unique morphism $e$ is also a total weak equivalence. According to what we have showed in \ref{prove C3}, $e = h \circ \lambda_{t, W}$ with $e_\triangleright = h$, where $h$ is the unique morphism such that the following diagram in $\calT_{B'\sqcup_{A'}C'}$ commutes:
\begin{equation*} \label{WEAK 2}
\tag{\ref{prove W2}.2}
\begin{gathered}
\xymatrix@R=1cm@C=0.1cm{
& \overline{u'}_!(v'_!(p(j)_!(X)))
     \ar@{=}[d]                                                  \\
& t_!(\overline{u}_!(v_!(X))) = t_!(\overline{v}_!(u_!(X)))
     \ar[dl]_-{t_!(\overline{u}_!(g_\triangleright))}
     \ar@{>->}[rr]^-{t_!(\overline{v}_!(f_\triangleright))}
     \ar'[d]_-{\overline{u'}_!(v'_!(j_\triangleright))}^-{\sim}[dd]
&  & t_!(\overline{v}_!(Y)) = \overline{v'}_!(p(k)_!(Y))
     \ar[dd]^(0.33){\overline{v'}_!(k_\triangleright)}_(0.33){\sim}
     \ar[dl]_-{t_!(a)}                              \\
  \overline{u'}_!(p(i)_!(Z)) = t_!(\overline{u}_!(Z))
     \ar@{>->}[rr]^(0.75){t_!(b)}
     \ar[dd]_(0.33){\overline{u'}_!(i_\triangleright)}^(0.33){\sim}
&  & t_!(W)
     \ar[dd]_(0.33){h}                                  \\
& \overline{u'}_!(v'_!(X')) = \overline{v'}_!(u'_!(X'))
     \ar[dl]_-{\overline{u'}_!(g'_\triangleright)}
     \ar@{>->}'[r][rr]^-{\overline{v'}_!(f'_\triangleright)}
&  & \overline{v'}_!(Y')
     \ar[dl]_-{a'}                          \\
  \overline{u'}_!(Z')
     \ar@{>->}[rr]^(0.75){b'}
&  & W'       }
\end{gathered}
\end{equation*}
Thus it suffices to prove that $t$ is a weak equivalence in $\calB$ and $h$ is a weak equivalence in $\calT_{B'\sqcup_{A'}C'}$.

Consider firstly the right commutative diagram in (\ref{WEAK 1}) in which both $u$ and $u'$ are cofibrations in $\calB$ and both the top and bottom squares are pushouts in $\calB$. Since $i$, $j$ and $k$ are total weak equivalences in $\calT$, it follows that $p(i)$, $p(j)$ and $p(k)$ are weak equivalences in $\calB$. But $\calB$ satisfies the axiom (W2), so $t$ is a weak equivalence in $\calB$, as desired.

Next, consider the commutative diagram (\ref{WEAK 2}), where both the top and bottom squares are pushouts in $\calT_{B'\sqcup_{A'}C'}$. Again, since $i$, $j$ and $k$ are total weak equivalences in $\calT$, it follows that $i_\triangleright$ is a weak equivalence in $\calT_{Z'}$, $j_\triangleright$ is a weak equivalence in $\calT_{X'}$, and $k_\triangleright$ is a weak equivalence in $\calT_{Y'}$. Thus $\overline{u'}_!(i_\triangleright)$, $\overline{u'}_!(v'_!(j_\triangleright))$ and
$\overline{v'}_!(k_\triangleright)$ are weak equivalences in $\calT_{B'\sqcup_{A'}C'}$, as the reindexing functors $\overline{u'}_!$, $v'_!$ and $\overline{v'}_!$ are exact. Since $\calT_{B'\sqcup_{A'}C'}$ satisfies the axiom (W2) as well, we conclude that $h$ is a weak equivalence in $\calT_{B'\sqcup_{A'}C'}$, as desired. \qed
\end{bfhpg}

\subsection{Two classical examples}
In this subsection we revist two classical examples given in \cite{Wald83} (see Example \ref{classical exa mor}), where the second one will be applied further in the proof of Theorem \ref{MAIN FOR REP CAT}.

For a Waldhausen category $\sfE$, recall that the morphism category $\Mor(\sfE)$ is also a Waldhausen category via a natural way; see Example \ref{classical exa mor}(1). Applying Theorem \ref{total E IS WALD} we can reobtain this Waldhausen structure via the codomain functor.

\begin{exa}\label{reobtain classical exa1}
Let $\sfE$ be a category.
Consider the codomain functor
\[
p: \Mor(\sfE) \longrightarrow \sfE, \ \ \mathrm{via}\ \ (f : X \to A) \longmapsto A.
\]
The following facts can be easily deduced:
\begin{itemize}

\item For each object $A$ in $\sfE$, the fiber $\Mor(\sfE)_A$ of $A$ is the overcategory $\sfE/A$.

\item The codomain functor $p$ is a Grothendieck opfibration. Indeed, for any object $(f: X \to A)$ in $\Mor(\sfE)$ and any morphism $u : A \to B$ in $\sfE$,
set $u_!((f: X \to A)) = (u \circ f : X \to B)$ and $\lambda_{u, f} = (\id_X, u)$; diagrammatically:
\[
\xymatrix@R=1cm@C=1cm{
   X
\ar@{=}[r]
\ar[d]_-{f}
&  X
\ar[d]^-{u \circ f}
\\
A
\ar[r]^-{u}
&  B}
\]
It is easy to check that $(\id_X, u)$ is a cocartesian morphism in $\Mor(\sfE)$ that is above $u$.

\item For any morphism $(a, u) : (f: X \to A) \to (g: Y \to B)$ in $\Mor(\sfE)$ that is above $u$, one has $(a, u)_\triangleright=(a, \id_B)$; diagrammatically:
\[
\xymatrix@R=0.5cm@C=0.35cm{
&&&&&
Y
\ar[dd]^(0.45){g}
\\
X
\ar@{=}[rrr]
\ar@/^1pc/[rrrrru]^(0.33){a}
\ar[dd]_(0.45){f}
&&&
X
\ar[rru]_{a}
\ar[dd]_(0.35){u \circ f}
\\
&&& \ar@/^0.2pc/[rr] &&
B
\\
A\ar[rrr]_{u}
\ar@/^0.3pc/@{}[rrru]^(0.6){u}
&&&
B
\ar@{=}[rru]}
\xymatrix@R=0.5cm@C=0.35cm{
\\ \quad \overset{p} \longmapsto \quad}
\xymatrix@R=0.5cm@C=0.35cm{
\\
&&&&&
B
\\
A\ar[rrr]_{u}
\ar@/^1pc/[rrrrru]^(0.33){u}
&&&
B
\ar@{=}[rru]}
\]

\item For any morphism $u: A \to B$ in $\sfE$ and any morphism $(a, \id_A): (f: X \to A) \to (f' : X' \to A)$ in $\sfE/A$, we have $u_!((a, \id_A)) = (a, \id_B)$ in $\sfE/B$.
\end{itemize}

Suppose now that $\sfE$ is a Waldhausen category with initial object $0$. Then from Example \ref{over cate is}(1), we see that for each object $A$ in $\sfE$, the overcategory $\sfE/ A$ also forms a Waldhausen category with initial object $(0 \to A)$, where a morphism $(a, \id_A)$ in $\sfE/ A$ is a cofibration (resp., weak equivalence) if $a$ is a cofibration (resp., weak equivalence) in $\sfE$. For any morphism $u : A \to B$ in $\sfE$, it is clear that $u_!$ preserves the initial objects, cofibrations, weak equivalences, and the pushout diagrams of axiom (C3). This implies that $u_!$ is exact from $\sfE/A$ to $\sfE/B$. Thus, $p$ is a Waldhausen opfibration.

The above facts ensure that we can apply Theorem \ref{total E IS WALD} to guarantee that $\Mor(\sfE)$ forms a Waldhausen category, where a morphism $(a, u) : (f: X \to A) \to (g: Y \to B)$ in $\Mor(\sfE)$ is a total cofibration (resp., total weak equivalence) if $u$ is a cofibration (resp., weak equivalence) in $\sfE$ and
$(a, u)_\triangleright=(a, \id_B)$ is a cofibration (resp., weak equivalence) in $\sfE/B$, that is, both $u$ and $a$ are cofibrations (resp., weak equivalences) in $\sfE$.

We mention that the above Waldhausen structure on $\Mor(\sfE)$ is identical with
the classical one; see Example \ref{classical exa mor}(1).
\qed
\end{exa}

Let $\sfE$ be a Waldhausen category. Recall that $\cMor(\sfE)$ is the full subcategory of $\Mor(\sfE)$ whose objects are cofibrations in $\sfE$. By \cite{Wald83}, $\cMor(\sfE)$ is a Waldhausen category; see Example \ref{classical exa mor}(2) for details. In the following example, applying Theorem \ref{total E IS WALD} again,
we reobtain this Waldhausen structure on $\cMor(\sfE)$ via the domain functor. Recall from Example \ref{over cate is}(2) that for any object $A$ in $\sfE$,
$\overline{A/\sfE}$ denotes the full subcategory of the undercategory $A/\sfE$ consisting of cofibrations $(A \rightarrowtail X)$ in $\sfE$.

\begin{exa}\label{rebtain classical 2}
Let $\sfE$ be a Waldhausen category with initial object $0$. Consider the domain functor 
\[
p: \cMor(\sfE) \longrightarrow \sfE, \ \ \mathrm{via}\ \
(f: A \rightarrowtail X) \longmapsto A.
\]
Note that $p$ is surjective on objects
as $\id_A \in \cMor(\sfE)$ for each object $A$ in $\sfE$.
We can deduce the following facts:

\begin{itemize}
\item For each object $A$ in $\sfE$,
the fiber $\cMor(\sfE)_A$ of $A$ is the category $\overline{A/\sfE}$.

\item The domain functor $p$ is a Grothendieck opfibration. Indeed, for any object $(f : A \rightarrowtail X)$ in $\cMor(\sfE)$ and any morphism $u: A \to B$ in $\sfE$, by the axiom (C3), there exists a pushout
\[\xymatrix@R=1cm@C=1cm{
  A
  \ar@{>->}[d]_-{f}
  \ar@{}[rd]^(0.6)>>{\lrcorner}
  \ar[r]^-{u}
& B
  \ar@{>->}[d]^-{\overline{f}}   \\
  X
  \ar[r]^-{\overline{u}}
& B \sqcup_A X                        }
\]
in $\sfE$ with $\overline{f}$ a cofibration, that is, $(\overline{f} : B \rightarrowtail B \sqcup_A X)$ is also an object in $\cMor(\sfE)$. Set $u_!((f: A \rightarrowtail X)) = (\overline{f} : B \rightarrowtail B +_A X)$ and $\lambda_{u, f} = (u, \overline{u})$. Then for each morphism $v : B \to C$ in $\sfE$ and each morphism $(v \circ u, c) : (f: A \rightarrowtail X) \to (h: C \rightarrowtail Z)$ in $\cMor(\sfE)$, by the universal property of the pushout, there exists a unique morphism $B +_A X \to Z$ such that the left diagram in $\sfE$ commutes:
\[
\xymatrix@R=0.5cm@C=0.5cm{
&&&&&
C
\ar@{>->}[dd]^(0.25){h}
\\
A
\ar[rrr]^-{u}
\ar@{}[rrrdd]|>>{\lrcorner}
\ar@/^1pc/[rrrrru]^(0.43){v \circ u}
\ar@{>->}[dd]_(0.25){f}
&&&
B
\ar[rru]^-{v}
\ar@{>->}[dd]^(0.25){\overline{f}}
\\
&&& \ar@/^0.2pc/[rr] &&
Z
\\
X
\ar[rrr]^(0.53){\overline{u}}
\ar@/^0.35pc/@{}[rrru]^(0.75){c}
&  &&
B \sqcup_A X
\ar@{.>}[rru]}
\xymatrix@R=0.5cm@C=0.35cm{
\\ \quad \overset{p} \longmapsto \quad}
\xymatrix@R=0.5cm@C=0.5cm{
\\
&&&&&
C
\\
A\ar[rrr]^-{u}
\ar@/^1pc/[rrrrru]^(0.43){v \circ u}
&&&
B
\ar[rru]^{v} }
\]
Thus, $(u, \overline{u})$ is a cocartesian morphism in $\cMor(\sfE)$ that is above $u$.

\item For any morphism $(u, a): ( f: A \rightarrowtail X) \to (g: B \rightarrowtail Y)$ in $\cMor(\sfE)$ that is above $u$, the fiber morphism $(u, a)_\triangleright=(\id_B, h)$, where $h$ is the unique morphism such that the diagram
    \begin{equation*} \label{pushout induce 1}
\tag{\ref{rebtain classical 2}.1}
\begin{gathered}
\xymatrix@R=0.75cm@C=0.5cm{
A
\ar@{>->}[dd]_-{f}
\ar[rrrr]^-{u}
&&&&
B
\ar@{>->}[dll]_(0.55){\overline{f}}
\ar@{>->}[dd]^-{g}                            \\
&& B \sqcup_A X
\ar@{.>}[drr]^(0.55){h}           \\
X
\ar[urr]^-{\overline{u}}
\ar[rrrr]^-{a}
&&&& Y
}
\end{gathered}
\end{equation*}
commutes.
\end{itemize}

By Example \ref{over cate is}(2), for any object $A$ in $\sfE$, the category $\overline{A/\sfE}$ forms a Waldhausen category with initial object $(\id_A : A \rightarrowtail A)$, where a morphism $(\id_A, a)$ in $\overline{A/\sfE}$ is a cofibration (resp., weak equivalence) if $a$ is a cofibration (resp., weak equivalence) in $\sfE$. For any morphism $u : A \to B$ in $\sfE$, we conclude that $u_!$ preserves the initial objects. For any morphism $(\id_A, g): (f: A \rightarrowtail X) \to (f' : A \rightarrowtail X')$ in $\overline{A/\sfE}$, consider the following left commutative diagram of solid arrows in $\sfE$, where both the front and back squares are pushouts:
\begin{equation*} \label{pushout induce 2}
\tag{\ref{rebtain classical 2}.2}
\begin{gathered}
\xymatrix@R=0.75cm@C=1cm{
& A
     \ar@{}[rrdd]^(0.6)>>{\lrcorner}
     \ar[rr]^(0.65){u}
     \ar@{>->}'[d]^-{f'}[dd]
&  & B
     \ar@{>->}[dd]^(0.33){\overline{f'}}                            \\
A
     \ar@{}[rrdd]^(0.6)>>{\lrcorner}
     \ar@{=}[ur]
     \ar[rr]^(0.66){u}
     \ar@{>->}[dd]_(0.33){f}
&  & B
     \ar@{=}[ur]
     \ar@{>->}[dd]^(0.33){\overline{f}}                                  \\
& X'
     \ar@{>->}'[r][rr]^-{\overline{u}'}
&  & B \sqcup_A X'                                                        \\
X
     \ar[ur]^-{g}
     \ar@{>->}[rr]_(0.66){\overline{u}}
&  & B \sqcup_A X
     \ar@{.>}[ur]_-{g'}     }
\xymatrix@R=1.25cm@C=1cm{ \\ \overset{p}\longmapsto }
\xymatrix@R=0.75cm@C=1cm{\\
& A
     \ar[rr]^(0.65){u}
&  & B                              \\
A
     \ar@{=}[ur]
     \ar[rr]^(0.66){u}
&  & B
     \ar@{=}[ur]                                 }
\end{gathered}
\end{equation*}
By the universal property of the front pushout, there exists a unique morphism $g': B \sqcup_A X \to B \sqcup_A X'$ such that the left diagram commutes, so we have $u_!((\id_A, g)) =(\id_B, g')$. By the pasting lemma, the bottom square in (\ref{pushout induce 2}) is also a pushout. Thus one can easily check that $u_!$ preserves the pushout diagrams of axiom (C3). Moreover, since $\sfE$ satisfies axiom (C3), it follows that if $(\id_A, g)$ is a cofibration in $\overline{A/\sfE}$ (that is, $g$ is a cofibration in $\sfE$), then so is $g'$, and hence $(\id_B, g')$ is a cofibration in $\overline{B/\sfE}$, which implies that $u_!$ preserves cofibrations.
If $(\id_A, g)$ is a weak equivalence in $\overline{A/\sfE}$ (that is, $g$ is a weak equivalence in $\sfE$), then by the axiom (W2), so is $g'$, and hence $(\id_B, g')$ is a weak equivalence in $\overline{B/\sfE}$. This implies that $u_!$ preserves weak equivalences as well. Thus, $u_!$ is exact from $\overline{A/\sfE}$ to $\overline{B/\sfE}$. According to what we showed above, we see that $p$ is Waldhausen opfibration.

The above facts allow us to apply Theorem \ref{total E IS WALD} to guarantee that $\cMor(\sfE)$ forms a Waldhausen category, where a morphism $(u, a) : (f : A \rightarrowtail X) \to (g : B \rightarrowtail Y)$ in $\cMor(\sfE)$ is a total cofibration (resp., total weak equivalence) if $u$ is a cofibration (resp., weak equivalence) in $\sfE$ and $(u, a)_\triangleright=(\id_B, h)$ appearing in (\ref{pushout induce 1}) is a cofibration (resp., weak equivalence) in $\overline{B/\sfE}$,
that is, both $u$ and $h$ are cofibrations (resp., weak equivalences) in $\sfE$.

We mention that the above Waldhausen structure on $\cMor(\sfE)$ is identical with the classical one; see Example \ref{classical exa mor}(2).
\qed
\end{exa}

\section{Proof of Theorem \ref{MAIN FOR REP CAT}}
\label{Appl in rep cat}
\noindent
We give the proof of Theorem \ref{MAIN FOR REP CAT} in this section as well as an example in the context of model categories to illustrate the application of this theorem. We begin by recalling the definition of left rooted quivers.

\subsection*{Left rooted quivers}\label{key prp for lrq}
Let $Q$ be a quiver with vertex set $V$ and arrow set $\Gamma$. For an arrow $\alpha$,  write $s(\alpha)$ for its source and $t(\alpha)$ for its target. By Enochs, Oyonarte and Torrecillas \cite{EOT04}, there exists a transfinite sequence $\{V_{\mu}\}_{\mu \, \mathrm{ordinal}}$ of subsets of $V$ as follows:
\begin{itemize}
\item for the first ordinal $\mu=0$, set $V_0=\emptyset$;

\item for a successor ordinal $\mu+1$, set
\[
V_{\mu+1}=\{i\in V\,|\,i\ \text{is not the target of any arrow } \alpha \ \text{with}\ s(\alpha)\notin\cup_{\gamma\leq\mu}V_{\gamma}\};
\]

\item for a limit ordinal $\mu$, set $V_{\mu}=\cup_{\gamma<\mu}V_{\gamma}$.
\end{itemize}

\begin{rmk}\label{main pro for ts}
It is clear that $V_1$ consists of all vertexes $i$ such that there is no arrow $\alpha$ with $t(\alpha)=i$. By \lemcite[2.7]{HJ19} and \corcite[2.8]{HJ19}, there is a chain
\[
V_0 \subseteq V_1 \subseteq \cdots \subseteq V_{\mu} \subseteq V_{\mu+1} \subseteq \cdots \subseteq V,
\]
and if $\alpha: i\to j$ is an arrow with $j\in V_{\mu+1}$ for some ordinal $\mu$, then $i$ must be in $V_{\mu}$.
\end{rmk}

\begin{dfn}
Let $Q = (V, \Gamma)$ be a quiver and $\{V_{\mu}\}_{\mu \, \mathrm{ordinal}}$ the transfinite sequence of subsets of $V$. Then $Q$ is said to be \emph{left rooted} \cite{EOT04} if there exists an ordinal $\zeta$ such that $V_{\zeta} = V$.
\end{dfn}

\begin{rmk}\label{exa for ts}
By \prpcite[3.6]{EOT04}, a quiver is left rooted if and only if it has no infinite sequence of arrows of the form
\[
\cdots \to \bullet\to \bullet \to \bullet
\]
where vertices are not necessarily different. It is easy to see that a left rooted quiver has no loops or oriented cycles. Moreover, by \cite[Remark 3.5 and Proposition 3.7]{DLLY2023}, the free category associated to a left rooted quiver is a direct category; see for instance Hovey \cite{Ho99}.
\end{rmk}

\begin{setup*}
Throughout this section,
\begin{itemize}
\item
let $Q = (V, \Gamma)$ be a left rooted quiver and $\{V_{\mu}\}_{\mu \leqslant \zeta}$ the transfinite sequence of subsets of $V$.

\item
let $Q_{\mu} = (V_{\mu}, \Gamma_{\mu})$ denote the subquiver of $Q$ spanned by $V_{\mu}$. In particular, $Q = Q_{\zeta}$. Note that every  $Q_{\mu}$ inherits a left rooted structure from $Q$.

\item
let $\sfE$ be a Waldhausen category with small coproducts, $\sfC$ (resp., $\sfW$) the class of all cofibrations (resp, weak equivalences) in $\sfE$, and suppose that $(\sfC, \sfC^\Box)$ forms a weak factorization system in $\sfE$.

\item let $\coE$ be the subcategory of $\sfE$ sharing the same objects whose morphisms are in $\sfC$, and
\[
\iota_\mu^* : \Rep{Q_{\mu+1}}{\coE} \to \Rep{Q_{\mu}}{\coE}
\]
be the restriction functor induced by the inclusion $\iota_\mu : Q_\mu \to Q_{\mu+1}$.
\end{itemize}
\end{setup*}

The proof of Theorem \ref{MAIN FOR REP CAT} goes through a transfinite induction and the most technical part of the proof is the induction step for successor ordinals. We will show that $\iota_\mu^*$ is a Waldhausen opfibration, and further, use Theorem \ref{total E IS WALD} to prove that if $\Rep{Q_\mu}{\coE}$ is a Waldhausen category via the construction given in Theorem \ref{MAIN FOR REP CAT}, then $\Rep{Q_{\mu+1}}{\coE}$ also forms a Waldhausen category via the same construction.

In order to prove that $\iota_\mu^*$ is a Waldhausen opfibration, we describe firstly the fiber $\Rep{Q_{\mu+1}}{\coE}_A$ of a representation $A$ in $\Rep{Q_{\mu}}{\coE}$.

\begin{ipg}\label{1-1 for exten}
Fix a representation $A$ in $\Rep{Q_{\mu}}{\coE}$ and let $X$ be a representation in $\Rep{Q_{\mu+1}}{\coE}_A$. For any vertex $i \in V_{\mu + 1} \backslash V_\mu$ and any arrow $\alpha \in \Gamma_{\mu +1}(\bullet, i)$, one has $X_{s(\alpha)}=\iota_\mu^*(X)_{s(\alpha)}=A_{s(\alpha)}$ as $s(\alpha) \in V_\mu$ by Remark \ref{main pro for ts}. By the universal property of the coproduct, there exists a unique morphism $\varphi^{X}_i$ such that the diagram
\begin{equation*} \label{commu 1}
\tag{\ref{1-1 for exten}.1}
\begin{gathered}
\xymatrix@R=1cm@C=2cm{
  X_{s(\alpha)} = A_{s(\alpha)}
   \ar[d]_-{\iota_\alpha}
   \ar@{>->}[dr]^-{X_\alpha}    \\
  \oplus_{\alpha \in \Gamma_{\mu +1}(\bullet, i)} A_{s(\alpha)}
   \ar[r]_-{\varphi^{X}_i}
& X_i
}
\end{gathered}
\end{equation*}
in $\sfE$ commutes, where $\iota_\alpha$ is the canonical injection. Hence, one gets a family of morphisms
\[
\{\varphi^{X}_i : \oplus_{\alpha \in \Gamma_{\mu +1}(\bullet, i)} A_{s(\alpha)}
                 \longrightarrow X_i \}_{i \in V_{\mu + 1} \backslash V_\mu}
\]
in $\sfE$.

We show next that $\varphi^{X}_i \in \sfC$ for each $i \in V_{\mu + 1} \backslash V_\mu$. It suffices to show that $\varphi^{X}_i$ has the left lifting property with respect to morphisms in $\sfC^\Box$ as $(\sfC, \sfC^\Box)$ is a weak factorization system by assumption. We explicitly construct a lift in the commutative diagram
\begin{equation*} \label{commu 2}
\tag{\ref{1-1 for exten}.2}
\begin{gathered}
\xymatrix@R=1cm@C=1cm{
  \oplus_{\alpha \in \Gamma_{\mu +1}(\bullet, i)} A_{s(\alpha)}
   \ar[d]_-{\varphi^{X}_i}
   \ar[r]^-{\sigma}
& M
   \ar[d]^-{\varepsilon}                          \\
X_i
   \ar[r]^-{\omega}
& N
}
\end{gathered}
\end{equation*}
in $\sfE$ with $\varepsilon \in \sfC^\Box$. For any arrow $\alpha \in \Gamma_{\mu +1}(\bullet, i)$, combining with (\ref{commu 1}), one has the following commutative diagram of solid arrows:
\[
\xymatrix@R=1cm@C=2cm{
  A_{s(\alpha)}
   \ar[d]_-{\iota_\alpha}
   \ar[r]^(0.4){\sigma \circ \iota_\alpha}
   \ar@{>->}@/_4pc/[dd]_-{X_\alpha}
& M
   \ar@{=}[d]                                        \\
  \oplus_{\alpha \in \Gamma_{\mu +1}(\bullet, i)} A_{s(\alpha)}
   \ar[d]_-{\varphi^{X}_i}
   \ar[r]^(0.45){\sigma}
& M
   \ar@{>>}[d]^-{\varepsilon}                          \\
  X_i
   \ar[r]^(0.4){\omega}
   \ar@{.>}[uur]|(0.66){\, \delta \,}
& N
}
\]
Since $X_\alpha \in \sfC$, there exists a lift $\delta$ in the outermost square, so $\varepsilon \circ \delta = \omega$ and $\delta \circ \varphi^{X}_i \circ \iota_\alpha = \sigma \circ \iota_\alpha$. By the universal property of the coproduct, one has $\delta \circ \varphi^{X}_i = \sigma$.
Thus, $\delta$ is the desired lift for (\ref{commu 2}).

Conversely, given a family $\{X_i\}_{i\in V_{\mu+1}\backslash V_\mu}$ of objects in $\sfE$, and a family
\[
\{\varphi_i: \oplus_{\alpha \in \Gamma_{\mu +1}(\bullet, i)} A_{s(\alpha)} \rightarrowtail X_i \}_{i \in V_{\mu + 1} \backslash V_\mu}
\]
of morphisms in $\sfC$, for any vertex $j \in V_{\mu + 1}$ and any arrow $\beta \in \Gamma_{\mu +1}$, one can define a representation $X$ in $\Rep{Q_{\mu+1}}{\sfE}$
as follows:
\begin{equation*}\label{dfn of widehat{X}}
\tag{\ref{1-1 for exten}.3}
X_j =
\begin{cases}
   A_j                   & \text{if } j \in V_\mu,      \\
   X_j                   & \text{if } j \in V_{\mu + 1} \backslash V_\mu;
\end{cases} \quad \quad \quad
X_\beta =
\begin{cases}
A_\beta                               & \text{if } t(\beta) \in V_\mu,      \\
\varphi_{t(\beta)} \circ \iota_\beta  & \text{if } t(\beta) \in V_{\mu + 1} \backslash V_\mu,
\end{cases}
\end{equation*}
where $\iota_\beta : A_{s(\beta)} \to \oplus_{\gamma \in \Gamma_{\mu +1}(\bullet, t(\beta))} A_{s(\gamma)}$ is the canonical injection. If $t(\beta) \in V_\mu$, then by Remark \ref{main pro for ts} one has $s(\beta) \in V_\mu$, so $\beta \in \Gamma_{\mu}$. In this case $X_\beta = A_\beta \in \sfC$. We claim that if $t(\beta) \in V_{\mu + 1} \backslash V_\mu$, then $X_\beta = \varphi_{t(\beta)} \circ \iota_\beta \in \sfC$, which implies that $X \in \Rep{Q_{\mu+1}}{\coE}$.

To prove the claim, it suffices to show that $\iota_\beta \in \sfC$ as $\sfC$ is closed under compositions. This is clear. Indeed, note that
\[
\xymatrix@R=1cm@C=1cm{
  0
   \ar[d]
   \ar[r]
& A_{s(\beta)}
   \ar[d]^-{\iota_\beta}                          \\
  \oplus_{\gamma \in \Gamma_{\mu +1}(\bullet, t(\beta))\backslash \{ \beta \}} A_{s(\gamma)}
   \ar[r]^-{\iota}
& \oplus_{\gamma \in \Gamma_{\mu +1}(\bullet, t(\beta))} A_{s(\gamma)}
}
\]
is a pushout in $\sfE$, where $0$ is the initial object in $\sfE$ and $\iota$ is the canonical injection. By axiom (C2), $0 \to \oplus_{\gamma \in \Gamma_{\mu +1}(\bullet, t(\beta))\backslash \{ \beta \}} A_{s(\gamma)}$ is contained in $\sfC$, so $\iota_\beta$ is contained in $\sfC$ as well by axiom (C3).

Moreover, it is clear that $X$ restricts to $A$. Thus, $X \in \Rep{Q_{\mu+1}}{\coE}_A$.
\end{ipg}

The above argument implies the following result:

\begin{lem}\label{1-1 for exten lem}
Let $A$ be a representation in $\Rep{Q_{\mu}}{\coE}$. Then there exists a bijective correspondence between representations in $\Rep{Q_{\mu+1}}{\coE}_A$ and families
\[
\{\oplus_{\alpha \in \Gamma_{\mu +1}(\bullet, i)} A_{s(\alpha)} \rightarrowtail \bullet \}_{i \in V_{\mu + 1}\backslash V_\mu}
\]
of morphisms in $\sfC$.
\end{lem}

\begin{rmk}\label{simple notation}
In a direct category, there is an induction procedure controlled by latching objects and functors; see for instance \cite{Ho99}. Note that the free category associated to $Q_{\mu}$ is a direct category. For any representation $A$ in $\Rep{Q_{\mu}}{\coE}$ and any vertex $i \in V_{\mu + 1} \backslash V_\mu$,
we emphasize that the object $\oplus_{\alpha \in \Gamma_{\mu +1}(\bullet, i)} A_{s(\alpha)}$ is nothing but the latching object of $A$ at $i$.
\end{rmk}

\begin{nota*}
Let $i$ be a vertex in $V_{\mu + 1} \backslash V_\mu$. To simplify the notation, from now on,
\begin{itemize}
\item
denote $\oplus_{\alpha \in \Gamma_{\mu +1}(\bullet, i)} A_{s(\alpha)}$ by $\sfL_i(A)$ for a representation $A \in \Rep{Q_{\mu}}{\coE}$;

\item
denote $\oplus_{\alpha \in \Gamma_{\mu +1}(\bullet, i)} u_{s(\alpha)}$ by $\sfL_i(u)$ for a morphism $u : A \to B \in \Rep{Q_{\mu}}{\coE}$.
\end{itemize}
\end{nota*}

For each representation $A$ in $\Rep{Q_{\mu}}{\coE}$ and each vertex $i$ in $V_{\mu+1} \backslash V_{\mu}$, recall that $\overline{\sfL_i(A)/\sfE}$ is the full subcategory of $\sfL_i(A)/\sfE$ consisting of cofibrations $(\sfL_i(A) \rightarrowtail \bullet)$ in $\sfE$; it is a Waldhausen category via the construction described in Example \ref{over cate is}(2). Define a rule
\[
\scrR : \Rep{Q_{\mu+1}}{\coE}_A \longrightarrow
\prod_{i \in V_{\mu+1} \backslash V_{\mu}} \overline{\sfL_i(A)/\sfE}
\]
as follows:
\begin{itemize}
\item For a representation $X$ in $\Rep{Q_{\mu+1}}{\coE}_A$, set $\scrR(X)$ to be the family
      \[ \{\varphi^{X}_i : \sfL_i(A) \rightarrowtail X_i \}_{i \in V_{\mu + 1} \backslash V_{\mu}}\]
      of morphisms in $\sfC$, where $\varphi^{X}_i$ is given in (\ref{commu 1}).

\item Let $f: X \to X'$ be a morphism in $\Rep{Q_{\mu+1}}{\coE}_A$. Then for any vertex $i \in V_{\mu + 1} \backslash V_{\mu}$ and any arrow $\alpha \in \Gamma_{\mu + 1}(\bullet, i)$, one has
              $$f_i \circ \varphi^{X}_i \circ \iota_\alpha
             = f_i \circ X_\alpha
             = X'_\alpha
             = \varphi^{X'}_i \circ \iota_\alpha,$$
 where the first and last equalities hold by (\ref{commu 1}), and  the second equality holds as $f_{s(\alpha)}=\iota_\mu^*(f)_{s(\alpha)}$ is the identity.
By the universal property of the coproduct, one has $f_i \circ \varphi^{X}_i = \varphi^{X'}_i$, that is, the bottom triangle in the following diagram commutes:
\[
\xymatrix@R=0.5cm@C=1.5cm{
   X_{s(\alpha)}
     \ar@{=}[r]
     \ar@{>->}[dd]_-{X_\alpha}
&  A_{s(\alpha)}
     \ar@{=}[r]
     \ar[d]_-{\iota_\alpha}
&  X'_{s(\alpha)}
     \ar@{>->}[dd]^-{X'_\alpha}             \\
&  \sfL_i(A)
     \ar@{>->}[dl]_-{\varphi^{X}_i}
     \ar@{>->}[dr]^-{\varphi^{X'}_i}         \\
   X_i
     \ar[rr]^-{f_i}
&& X'_i
}
\]
Consequently, $f_i$ is a morphism in $\overline{\sfL_i(A)/\sfE}$, and we set $\scrR(f)$ to be the family $\{f_i\}_{i \in V_{\mu+1} \backslash V_{\mu}}$.
\end{itemize}

It is routine to check that $\scrR$ is a functor. Furthermore, one can easily deduce the following result from Lemma \ref{1-1 for exten lem}.

\begin{prp}\label{isomorphism of fibe}
The functor $\scrR: \Rep{Q_{\mu+1}}{\coE}_A \to
\prod_{i \in V_{\mu+1} \backslash V_{\mu}} \overline{\sfL_i(A)/\sfE}$ given above is an isomorphism.
\end{prp}

With help of Lemma \ref{1-1 for exten lem}, we get the next result,
which exhibits the opfibrational nature of $\iota_\mu^*$.

\begin{lem}\label{iota is Gro op}
The restriction functor $\iota_\mu^*: \Rep{Q_{\mu+1}}{\coE} \to \Rep{Q_{\mu}}{\coE}$ is a Grothendieck opfibration.
\end{lem}

\begin{prf*}
Let $u: A \to B$ be a morphism in $\Rep{Q_{\mu}}{\coE}$ and $X$ a representation in $\Rep{Q_{\mu+1}}{\coE}_A$. We construct next a representation $u_!(X)$ in $\Rep{Q_{\mu+1}}{\coE}_B$ and a cocartesian morphism $\lambda_{u, X}: X \to u_!(X)$ in $\Rep{Q_{\mu+1}}{\coE}$ which is above $u$. Indeed, from \ref{1-1 for exten}, we see that $X$ induces a family $\{\varphi^{X}_i : \sfL_i(A) \rightarrowtail X_i \}_{i \in V_{\mu + 1} \backslash V_{\mu}}$ of morphisms in $\sfC$. For each vertex $i \in V_{\mu + 1} \backslash V_\mu$, since $\sfE$ satisfies the axiom (C3), it follows that there exists a pushout
\begin{equation*} \label{pushout d 1}
\tag{\ref{iota is Gro op}.1}
\begin{gathered}
\xymatrix@R=1cm@C=1cm{
     \sfL_i(A)
       \ar[r]^-{\sfL_i(u)}
       \ar@{>->}[d]_-{\varphi^{X}_i}
&    \sfL_i(B)
       \ar@{>->}[d]^-{\varphi_i}             \\
     X_i
        \ar[r]^-{\theta_i}
&    \sfL_i(B) \sqcup_{\sfL_i(A)} X_i
}
\end{gathered}
\end{equation*}
in $\sfE$ with $\varphi_i \in \sfC$, from which one obtains a family $\{\varphi_i \}_{i \in V_{\mu + 1} \backslash V_\mu}$ of morphisms in $\sfC$. By \ref{1-1 for exten} again, one gets a representation $u_!(X) \in \Rep{Q_{\mu+1}}{\coE}_B$ (see (\ref{dfn of widehat{X}}) for its definition). For any vertex $j \in V_{\mu + 1}$, set
\begin{equation*}
({\lambda_{u, X}})_j =
\begin{cases}
u_j         & \text{if } j \in V_{\mu};      \\
\theta_j    & \text{if } j \in V_{\mu + 1} \backslash V_\mu.
\end{cases}
\end{equation*}
It is routine to check that $\lambda_{u, X}$ is a natural transformation from $X$ to $u_!(X)$, or equivalently, a morphism in $\Rep{Q_{\mu+1}}{\coE}$. Note that for any $j\in V_{\mu}$, one has $(\iota_\mu^*({\lambda_{u, X}}))_j=({\lambda_{u, X}})_j=u_j$, so $\lambda_{u, X}$ is above $u$.

Next, we prove that $\lambda_{u, X}$ is cocartesian. Consider the following diagram of solid arrows
\[
\xymatrix@R=1cm@C=2cm{
&&  Z  \\
    X
      \ar[r]_-{\lambda_{u, X}}
      \ar@/^1pc/[rru]^(0.43){g}
&   u_!(X)
      \ar@{.>}[ru]_-{h}
}
\xymatrix@R=0.5cm@C=0.35cm{
\\ \quad \overset{\iota_\mu^*} \longmapsto \quad}
\xymatrix@R=1cm@C=2cm{
&&  C  \\
    A
      \ar[r]_-{u}
      \ar@/^1pc/[rru]^(0.43){v \circ u}
&   B
      \ar[ru]_-{v}
}
\]
in which $g$ is above $v \circ u$. We need to construct a unique morphism $h$ above $v$ such that the left triangle in $\Rep{Q_{\mu+1}}{\coE}$ commutes.  Since (\ref{pushout d 1}) is a pushout for each vertex $i \in V_{\mu + 1} \backslash V_\mu$, it follows that there exists a unique morphism $\delta_i$ such that
the diagram
\[
\xymatrix@R=0.5cm@C=0.5cm{
&&&&&
\sfL_i(C)
   \ar@{>->}[dd]^(0.6){\varphi^{Z}_i}
\\
\sfL_i(A)
   \ar[rrr]_(0.45){\sfL_i(u)}
   \ar@/^1pc/[rrrrru]^(0.66){\sfL_i(v \circ u)}
   \ar@{>->}[dd]_(0.7){\varphi^{X}_i}
&&&
\sfL_i(B)
   \ar[rru]_{\sfL_i(v)}
   \ar@{>->}[dd]^(0.7){\varphi_i}
\\
&&& \ar@/^0.25pc/^(0.25){g_i}[rr]
&&
Z_i
\\
X_i
   \ar[rrr]_(0.4){\theta_i}
   \ar@/^0.45pc/@{}[rrru]
&&&
\sfL_i(B) \sqcup_{\sfL_i(A)} X_i
    \ar@{.}[rru]_-{\delta_i}
}
\]
in $\sfE$ commutes. For any vertex $j \in V_{\mu + 1}$, set
\begin{equation*}
h_j =
\begin{cases}
v_j      & \text{if } j \in V_{\mu};      \\
\delta_j & \text{if } j \in V_{\mu + 1} \backslash V_\mu.
\end{cases}
\end{equation*}
It is easy to check that $h$ is above $v$ and $h \circ \lambda_{u, X} = g$.
\end{prf*}

\begin{rmk}\label{EXPLIAN FOR FIBER}
Let $f: X \to Y$ be a morphism in $\Rep{Q_{\mu+1}}{\coE}$ that is above a morphism $u: A \to B$ in  $\Rep{Q_\mu}{\coE}$. Under the identification through $\scrR$,
one may view $Y$ as the family
\[
\{\varphi^{Y}_i : \sfL_i(B)\rightarrowtail Y_i \}_{i \in V_{\mu + 1} \backslash V_\mu}
\]
of morphisms in $\sfC$ by Lemma \ref{1-1 for exten lem}. For any vertex $i \in V_{\mu + 1} \backslash V_\mu$, by the universal property of the pushout, there is a unique morphism $\rho_i$ such that the following diagram in $\sfE$ commutes:
\[
\xymatrix@R=0.4cm@C=0.5cm{
&&&&&
\sfL_i(B)
   \ar@{>->}[dd]^(0.6){\varphi^{Y}_i}
\\
\sfL_i(A)
   \ar[rrr]_(0.45){\sfL_i(u)}
   \ar@/^1pc/[rrrrru]^(0.66){\sfL_i(u)}
   \ar@{>->}[dd]_(0.7){\varphi^{X}_i}
&&&
\sfL_i(B)
   \ar@{=}[rru]
   \ar@{>->}[dd]^(0.7){\varphi_i}
\\
&&& \ar@/^0.25pc/^(0.25){f_i}[rr]
&&
Y_i
\\
X_i
   \ar[rrr]_(0.4){\theta_i}
   \ar@/^0.45pc/@{}[rrru]
&&&
\sfL_i(B) \sqcup_{\sfL_i(A)} X_i
    \ar@{.}[rru]_-{\rho_i}
}
\]
Hence, the fiber morphism $f_\triangleright$, up to the identification through $\scrR$, is the family
\[
\{\rho_i : \sfL_i(B) \sqcup_{\sfL_i(A)} X_i \longrightarrow Y_i \}_{i \in V_{\mu+1} \backslash V_{\mu}}.
\]
\end{rmk}

Now we can show that $\iota_\mu^*$ is a Waldhausen opfibration.

\begin{prp}\label{iota mu is a Waldop}
The restriction functor
$\iota_\mu^* : \Rep{Q_{\mu+1}}{\coE} \to \Rep{Q_{\mu}}{\coE}$
is a Waldhausen opfibration.
\end{prp}

\begin{prf*}
The functor $\iota_\mu^*$ is a Grothendieck opfibration by Lemma \ref{iota is Gro op}. For each representation $A \in \Rep{Q_{\mu}}{\coE}$  and each vertex $i \in V_{\mu+1} \backslash V_{\mu}$, note that $\overline{\sfL_i(A)/\sfE}$ is a Waldhausen category by Example \ref{over cate is}(2). It follows from Proposition \ref{isomorphism of fibe} that the fiber $\Rep{Q_{\mu+1}}{\coE}_A$ is also a Waldhausen category.

For any morphism $u : A \to B \in \Rep{Q_{\mu}}{\coE}$, note that the associated reindexing functor $u_! : \Rep{Q_{\mu+1}}{\coE}_A \to \Rep{Q_{\mu+1}}{\coE}_B$, up to identification through $\scrR$ again, is the family $\prod_{i \in V_{\mu+1} \backslash V_{\mu}} \sfL_i(u)_!$, where $\sfL_i(u)_!: \overline{\sfL_i(A)/\sfE}\to\overline{\sfL_i(B)/\sfE}$ is the reindexing functor of $\sfL_i(u)$ as described in Example \ref{rebtain classical 2}.
As in Example \ref{rebtain classical 2}, one can show that each $\sfL_i(u)_!$ is exact, so are $\prod_{i \in V_{\mu+1} \backslash V_{\mu}} \sfL_i(u)_!$ and $u_!$.
\end{prf*}

Now we give a proof of Theorem \ref{MAIN FOR REP CAT}. Recall from the paragraph before Theorem \ref{MAIN FOR REP CAT} that for any morphism $f : X \to Y$ in $\Rep{Q_{\mu}}{\coE}$ with $\mu \leqslant \zeta$ and any vertex $i$ in $V_\mu$, by considering the following commutative diagram in $\sfE$, where the inner square is a pushout, we obtain a natural morphism $\rho_i$:
\begin{equation}\label{ast}
\tag{$\ast$}
\begin{gathered}
\xymatrix@R=1cm@C=1cm{
   \sfL_i(X)
     \ar@{>->}[d]_{\varphi^X_i}
     \ar[r]^-{\sfL_i(f)}
     \ar@{}[rd]^(0.6)>>{\lrcorner}
& \sfL_i(Y)
     \ar@{>->}[d]_{}
     \ar@/^0.8pc/[ddr]^{}                          \\
  X_i
     \ar[r]^-{\overline{f}_i}
     \ar@/_0.8pc/[drr]_{f_i}
& \sfL_i(Y) \sqcup_{\sfL_i(X)} X_i
     \ar@{.>}[dr]|-{\, \rho_i \,}                              \\
&& Y_i.}
\end{gathered}
\end{equation}

\begin{bfhpg}[Proof of Theorem \ref{MAIN FOR REP CAT}]\label{pf of thB}
It is clear that $\Rep{Q}{\coE}$ admits an initial object, the zero functor. We proceed by a transfinite induction. Let $\mu$ be an ordinal with $\mu \leqslant \zeta$. If $\mu = 1$, then for any morphism $f: X \to Y \in \Rep{Q_1}{\coE}$, one has $\sfL_i (X) = 0 = \sfL_i (Y)$ by Remark \ref{main pro for ts}. Thus the induced morphism
\[
\rho_i: \sfL_i(Y)\sqcup_{\sfL_i(X)}X_i \to Y_i
\]
in (\ref{ast}) is precisely $f_i : X_i \to Y_i$. In this case, one defines $f$ as a cofibration (resp., weak equivalence) in $\Rep{Q_1}{\coE}$ if $f_i$ is a cofibration (resp., weak equivalence) in $\sfE$ for each vertex $i \in V_1$. It is clear that $\Rep{Q_1}{\coE}$ forms a Waldhausen category.

Now we carry our the induction procedure. There are two cases. If $\mu$ is limit ordinal, then each vertex $i \in V_\mu$ lies in some $V_\lambda$ with $\lambda < \mu$. By the induction hypothesis, one checks routinely that $\Rep{Q_\mu}{\coE}$ forms a Waldhausen category via the construction described in the statement of Theorem \ref{MAIN FOR REP CAT}.

For the case of non-limit ordinals, suppose that $\Rep{Q_\mu}{\coE}$ is a Waldhausen category by the construction given in the statement of Theorem \ref{MAIN FOR REP CAT}. Since $\iota_\mu^*: \Rep{Q_{\mu+1}}{\coE} \to \Rep{Q_\mu}{\coE}$ is a Waldhausen opfibration by Proposition \ref{iota mu is a Waldop}, it follows from Theorem \ref{total E IS WALD} that $\Rep{Q_{\mu+1}}{\coE}$ is a Waldhausen category. Moreover, because the fiber morphism $f_\triangleright$, up to identification through $\scrR$, is the family
\[
\{\rho_i: \sfL_i(B) \sqcup_{\sfL_i(A)} X_i  \longrightarrow Y_i \}_{i \in V_{\mu+1} \backslash V_{\mu}}
\]
(see Remark \ref{EXPLIAN FOR FIBER}), $\sfL_i(A) = \sfL_i(X)$ and $\sfL_i(B) = \sfL_i(Y)$, it follows that a morphism $f: X \to Y$ that is above $u: A \to B$ in $\Rep{Q_\mu}{\coE}$ is a cofibration (resp., weak equivalence) if and only if $u$ is a cofibration (resp., weak equivalence) in $\Rep{Q_\mu}{\coE}$ and the induced morphism $\rho_i : \sfL_i(Y) \sqcup_{\sfL_i(X)} X_i \to Y_i$ lies in $\sfC$ (resp. $\sfW$) for each vertex $i \in V_{\mu+1} \backslash V_{\mu}$. However, by the induction hypothesis, $u$ is a cofibration (resp., weak equivalence) in $\Rep{Q_\mu}{\coE}$ if and only if $\rho_i$ lies in $\sfC$ (resp. $\sfW$) for each vertex $i \in V_\mu$. Consequently, $f$ is a cofibration (resp., weak equivalence) in $\Rep{Q_{\mu+1}}{\coE}$ if and only if $\rho_i$ lies in $\sfC$ (resp. $\sfW$) for all vertexes $i \in V_{\mu + 1}$. This completes the proof.
\qed
\end{bfhpg}


\begin{rmk}\label{pointed w}
If $\sfW$ is closed under small coproducts and pushouts along morphisms in $\sfC$, and satisfies the 2-out-of-3 property, then a morphism $f : X \to Y$ in $\Rep{Q}{\coE}$ is a weak equivalence (that is, $\rho_i \in \sfW$ for each vertex $i \in V$) if and only if $f_i : X_i \to Y_i \in \sfW$ for each vertex $i \in V$.
To show this assertion, consider the commutative diagram (\ref{ast}). For each vertex $i \in V$, since $\varphi^X_i : \sfL_i(X) \to X$ is contained in $\sfC$ (see \ref{1-1 for exten}) and $\sfW$ is closed under small coproducts, we deduce that $\sfL_i(f)$ is contained in $\sfW$, so is $\overline{f}_i$. Thus $\rho_i$ is contained in $\sfW$ if and only if so is $f_i$, as $\sfW$ satisfies the 2-out-of-3 property.
\end{rmk}

Let $\sfM$ be a Quillen model category with small coproducts and $(\sfC, \sfW, \sfF)$ its model structure. Since the free category associated to $Q$ is a direct category, it follows that $\Rep{Q}{\sfM}$ is a Quillen model category (see Hirschhorn \cite{Hir03}), whose model structure $(\sfC\sfR, \sfW\sfR, \sfF\sfR)$
are described as follows:
\begin{itemize}
\item a morphism $f: X \to Y \in \sfC\sfR$ if the induced morphism $\sfL_i(Y)\sqcup_{\sfL_i(X)}X_i  \to Y_i$ in (\ref{ast}) is in $\sfC$ for each vertex $i \in V$;

\item a morphism $f: X \to Y \in \sfW\sfR$ if $f_i: X_i \to Y_i \in \sfW$ for each vertex $i \in V$;

\item a morphism $f : X \to Y \in \sfF\sfR$ if $f_i : X_i \to Y_i \in \sfF$ for each vertex $i \in V$.
\end{itemize}
In what follows, denote by $\sfCof$ (resp., $\sfCof\sfR$) the full subcategory of $\sfM$ (resp., $\Rep{Q}{\sfM}$) consisting of cofibrant objects.

By Example \ref{exa model} and Theorem \ref{MAIN FOR REP CAT}, the full subcategories $\sfCof\sfR$ and $\Rep{Q}{\overline{\sfCof}}$ of $\Rep{Q}{\sfM}$ form Waldhausen categories, where $\overline{\sfCof}$ is the subcategory of $\sfM$ consisting of cofibrant objects and cofibrations between them. As an illustration of Theorem \ref{MAIN FOR REP CAT}, we show in the following example that these two Waldhausen categories are actually coincident.

\begin{exa}\label{fur illusta on model}
By Example \ref{exa model}, $\sfCof\sfR$ is a Waldhausen category whose cofibrations (resp., weak equivalences) are cofibrations (resp., weak equivalences) between cofibrant objects in $\Rep{Q}{\sfM}$. On the other hand, by Example \ref{exa model} and Theorem \ref{MAIN FOR REP CAT}, $\Rep{Q}{\overline{\sfCof}}$ is also a Waldhausen category, where a morphism $f : X \to Y $ is a cofibration $($resp., weak equivalence$)$ if the induced morphism $\rho_i: \sfL_i(Y)\sqcup_{\sfL_i(X)}X_i \to Y_i$ in (\ref{ast}) is a cofibration (resp., weak equivalence) between cofibrant objects in $\sfM$ for each vertex $i \in V$.

We show that the above two Waldhausen categories coincide. A key observation is the following fact. Let $X$ be a representation in $\Rep{Q}{\sfM}$. For each vertex $i \in V$ and each arrow $\alpha \in \Gamma(\bullet, i)$, the following statements are equivalent:
\begin{enumerate}
\item the morphism $\varphi^X_i : \sfL_i(X) \to X_i$ is a cofibration in $\sfM$;

\item the morphism $X_\alpha : X_{s(\alpha)} \to X_i$ is a cofibration in $\sfM$;

\item $X$ is a cofibrant object in $\Rep{Q}{\sfM}$.
\end{enumerate}

The equivalence between (1) and (2) can be proved by an argument similar to that in \ref{1-1 for exten}. To see the equivalence between (1) and (3), note that $X$ is a cofibrant object in $\Rep{Q}{\sfM}$ if and only if $0\sfR \to X$ is a cofibration, and if and only if the morphism $\varphi^X_i : \sfL_i(X) \to X_i$ in the commutative diagram
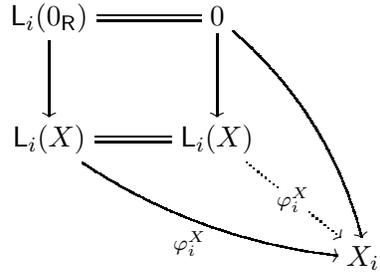
\[\xymatrix@R=1cm@C=1cm{
  \sfL_i(0\sfR)
     \ar[d]_{}
     \ar@{=}[r]^{}
& 0 \ar[d]_{}
     \ar@/^0.8pc/[ddr]^{}                          \\
  \sfL_i(X)
     \ar@{=}[r]^-{}
     \ar@/_0.8pc/[drr]_{\varphi^X_i}
& \sfL_i(X)
     \ar@{.>}[dr]|-{\, \varphi^X_i \,}                              \\
&& X_i.}\]
is a cofibration in $\sfM$.

Now we check the following facts:
\begin{itemize}
\item $\sfCof\sfR$ and $\Rep{Q}{\overline{\sfCof}}$ have the same objects. If $X$ is a representation in $\Rep{Q}{\overline{\sfCof}}$, then each $X_\alpha : X_{s(\alpha)} \to X_i$ is a cofibration in $\sfM$, so $X$ is contained in $\sfCof\sfR$ by the above equivalence. Conversely, given $X$ contained in $\sfCof\sfR$, then by the above equivalence, each $X_\alpha : X_{s(\alpha)} \to X_i$ is a cofibration in $\sfM$. Furthermore, since the functor $\sfL_i$ sends cofibrations in $\Rep{Q}{\sfM}$ to cofibrations in $\sfM$ by \corcite[5.1.5]{Ho99}, it follows that $0 \to \sfL_i(X)$ is a cofibration in $\sfM$, so is $0 \to X_i$, as $\sfC$ is closed under compositions. Thus each $X_i$ is a cofibrant object in $\sfM$. Consequently, $X$ is contained in $\Rep{Q}{\overline{\sfCof}}$.

\item It easily follows from the definitions that cofibrations in $\sfCof\sfR$ and $\Rep{Q}{\overline{\sfCof}}$ coincide.

\item Weak equivalences in $\sfCof\sfR$ and $\Rep{Q}{\overline{\sfCof}}$ are the same. By Remark \ref{pointed w}, it suffices to show that the class $\Mor(\sfCof) \cap \sfW$ is closed under small coproducts and pushouts along morphisms in $\sfC$, and satisfies the 2-out-of-3 property. Firstly, since $(\sfC \cap \sfW, \sfF)$ is a weak factorization system in $\sfM$, it follows that $\sfC \cap \sfW$ is closed under small coproducts by Remark \ref{need pro for wfs}(c). Thus the class of $\Mor(\sfCof) \cap \sfC \cap \sfW$ also has this property. By \corcite[7.7.2]{Hir03}, we conclude that $\Mor(\sfCof) \cap \sfW$ is also closed under small coproducts. Secondly, by \prpcite[2.29]{K-P97}, $\Mor(\sfCof) \cap \sfW$ is also closed under pushouts along cofibrations. Finally, it is clear that $\Mor(\sfCof) \cap \sfW$ satisfies the 2-out-of-3 property.

\end{itemize}
\end{exa}



\bibliographystyle{amsplain-nodash}


\def\cprime{$'$}
  \providecommand{\arxiv}[2][AC]{\mbox{\href{http://arxiv.org/abs/#2}{\sf
  arXiv:#2 [math.#1]}}}
  \providecommand{\oldarxiv}[2][AC]{\mbox{\href{http://arxiv.org/abs/math/#2}{\sf
  arXiv:math/#2
  [math.#1]}}}\providecommand{\MR}[1]{\mbox{\href{http://www.ams.org/mathscinet-getitem?mr=#1}{#1}}}
  \renewcommand{\MR}[1]{\mbox{\href{http://www.ams.org/mathscinet-getitem?mr=#1}{#1}}}
\providecommand{\bysame}{\leavevmode\hbox to3em{\hrulefill}\thinspace}
\providecommand{\MR}{\relax\ifhmode\unskip\space\fi MR }
\providecommand{\MRhref}[2]{%
  \href{http://www.ams.org/mathscinet-getitem?mr=#1}{#2}
}
\providecommand{\href}[2]{#2}

\end{document}